\documentclass[11pt]{amsart}
\usepackage{amssymb, amsmath, amsthm}
\usepackage[latin1]{inputenc}
\usepackage{graphicx}
\usepackage[final]{hyperref}
\usepackage[a4paper, centering]{geometry}
\usepackage{color}
\usepackage{graphicx} 

\newtheorem{theorem}{Theorem}
\newtheorem{proposition}[theorem]{Proposition}
\newtheorem{lemma}[theorem]{Lemma}

\theoremstyle{definition}
\newtheorem{definition}[theorem]{Definition}
\theoremstyle{remark}
\newtheorem{remark}[theorem]{Remark}

\def\R{\mathbb{R}}

\def\pscal#1#2{\left\langle#1,\,#2\right\rangle}
\def\dist{d_{\partial\Omega}}

\def\osubjet{J^{2, -}_{\Omega}}
\def\osuperjet{J^{2, +}_{\Omega}}

\def \e{\varepsilon}

\def\X{\mathbf{X}}
\def\Xe{\X_{\epsilon}}

\def\DN{\Delta^N_{\infty}}
\def\Dp{\Delta^+_{\infty}}
\def\Dm{\Delta^-_{\infty}}
\def\Dinf{\Delta_{\infty}}
\def\lmax{\lambda_{\text{max}}}
\def\lmin{\lambda_{\text{min}}}
\def\inradius{\rho_{\Omega}}
\def\PN{P_N}
\def\lN{\lambda_N} 
\def\uN{u_N}

\DeclareMathOperator{\Cut}{\overline \Sigma}  
\DeclareMathOperator{\high}{M}

\DeclareMathOperator{\argmax}{argmax}

\begin{document}

%%%%%%%%%%%%%%%%%%%%%%%%%%%%%%%%%%%%%%%%%%%%%%%%%%%%%%%%%%
%amsart format
\title[Boundary value problems for infinity Laplacian]%
{Characterization of stadium-like domains via boundary value problems for the infinity Laplacian}%
\author[G.~Crasta, I.~Fragal\`a]{Graziano Crasta,  Ilaria Fragal\`a}
\address[Graziano Crasta]{Dipartimento di Matematica ``G.\ Castelnuovo'', Univ.\ di Roma I\\
P.le A.\ Moro 2 -- 00185 Roma (Italy)}
\email{crasta@mat.uniroma1.it}

\address[Ilaria Fragal\`a]{
Dipartimento di Matematica, Politecnico\\
Piazza Leonardo da Vinci, 32 --20133 Milano (Italy)
}
\email{ilaria.fragala@polimi.it}

%49K20    Problems involving partial differential equations
%49K30    Optimal solutions belonging to restricted classes
%35J70 Degenerate elliptic equations
%35N25 Overdetermined boundary value problems
\date{September 5, 2015}

\begin{abstract} We give a complete characterization, as ``stadium-like  domains'',  of convex subsets $\Omega$ of $\mathbb{R}^n$ where a solution exists to Serrin-type overdetermined boundary value problems in which the operator is either the  infinity Laplacian or  its normalized version. In case of the not-normalized operator, our results extend those 
obtained in a previous work,  where the problem was solved under some geometrical restrictions on $\Omega$.
In case of the normalized operator, we also show that stadium-like domains are precisely the unique convex sets in $\mathbb{R}^n$ where the solution to a Dirichlet problem is of class $C^{1,1} (\Omega)$. 
\end{abstract}

\maketitle

\section{Introduction}
Consider the following 
Serrin-type problems for the infinity Laplace operator $\Dinf$ or its normalized version $\DN$: 
\begin{equation}\label{f:s}
\begin{cases} 
		-\Dinf u = 1 &\text{in}\ \Omega,\\
		u = 0 &\text{on}\ \partial\Omega, \\ 
		|\nabla u| = c &\text{on}\ \partial\Omega,
	\end{cases}
\end{equation}
and
\begin{equation}\label{f:sN}
	\begin{cases} 
		-\DN u = 1 &\text{in}\ \Omega,\\
		u = 0 &\text{on}\ \partial\Omega,\\ 
		|\nabla u| = c &\text{on}\ \partial\Omega\,. 
	\end{cases}
\end{equation}
Aim of this paper is to provide a complete characterization of convex domains $\Omega \subset \R ^n$ where such problems admit a solution.

Following the seminal paper by Serrin \cite{Se} and the huge amount of literature after it (see for instance \cite{BH, BrPr,  f, fg, fgk, GL, K2, Vog}), 
overdetermined boundary value problems involving the infinity Laplace operator were firstly considered
only few years ago by Buttazzo and Kawohl (see \cite{butkaw}). 
In fact, due to the high degeneracy of the operator, all the different methods exploited in the literature to obtain symmetry results for
overdetermined boundary value problems fail when applied to problems \eqref{f:s}-\eqref{f:sN}.

In \cite{butkaw}, Buttazzo and Kawohl dealt with a simplified version of problems  \eqref{f:s}-\eqref{f:sN}, 
which consists in looking for solutions having the same level lines as the distance function  to the boundary of $\Omega$, 
which are called {\it web--functions} (see Section~\ref{secprel} below). 
This simplification essentially reduces the problem to a one-dimensional setting, allowing to prove that the existence of a web--solution implies a precise geometric condition on $\Omega$, which is the coincidence of its {\it cut locus} and {\it high ridge} (see again Section~\ref{secprel} for the definitions). In particular,  such condition does not imply symmetry, at least if taken alone without any additional boundary regularity requirement. 
In our previous paper \cite{CFb} we studied the geometry of domains whose cut locus and high ridge agree, by providing a complete characterization of them in dimension $n=2$, and
in higher dimensions under convexity constraint; in particular, these results reveal that planar convex sets with the same cut locus and high ridge are tubular neighborhoods of a line segment (possibly degenerated into a point).  
Moreover, in \cite{CFc} we were able to carry over the study of problem \eqref{f:s} in the class of web--functions, by dropping all the regularity hypotheses on both the domain and the solution previously asked in \cite{butkaw}. 

The study of problems~\eqref{f:s}--\eqref{f:sN} in their full generality, namely without imposing the solution to be a  web-function, turns out to be much more challenging. As for problem \eqref{f:sN}, to our knowledge it has never been undergone. As for problem \eqref{f:s},  in a recent work we proved that, among convex sets, those having the same cut locus and high ridge - that we call ``stadium-like domains'' - are the only ones for which  {\it whatever} solution (not necessarily of web type) exists, see \cite[Thm.\ 5]{CFd}. 
As a drawback, we  needed to ask the following {\it a priori} geometrical hypothesis on the convex domain $\Omega$: there exists an inner ball,  of radius equal to the maximum of the distance from the boundary, touching $\partial \Omega$ at two diametral points. Moreover, we also needed the technical assumption that $\Omega$ satisfies an interior sphere condition at every point of the boundary. 

The approach we adopted for the proof relies on the study of a suitable $P$-function along the gradient flow of the unique solution to the Dirichlet problem. In particular, the diametral ball condition was used as a fundamental picklock to get the result. Indeed,  it allowed us to overcome the possible lack of regularity of the solution, which is  an intrinsic phenomenon; we refer to \cite[Sections 5 and 6]{CFd} for more details, including regularity thresholds. 

However, there was no reason to think that the geometric assumptions made on $\Omega$  
should be really necessary, so that the conclusion reached in \cite{CFd}
was not completely satisfactory. 

We can now introduce the contents of this paper, by describing its main results: 

\begin{itemize}
\item[--] Theorem \ref{t:s} improves the achievement of \cite[Thm.\ 5]{CFd}, by showing that it continues to hold  without any geometric assumptions on $\Omega$, i.e.\ when both the diametral ball assumption and the interior sphere condition are removed. Contrarily to our previous belief, it is possible to arrive at this conclusion by completely circumventing regularity matters, but rather exploiting the observation that a suitable web--function is always a super--solution to our problem (see Proposition \ref{l:superweb}). 

\smallskip
\item[--] Theorem \ref{t:sN} states that the same result
(in its  fully general version when no assumption is made on the convex set $\Omega$), 
holds true in the case of the normalized infinity Laplace operator $\DN$.
Recently, such operator has attracted an increasing interest for its applications and connections with different areas, in particular ``tug of war'' differential games \cite{ArmSm2012, KoSe, PSSW}. 
As we pointed out in \cite[Remark 3]{CFd}, in order to deal with problem \eqref{t:sN} a missing key ingredient was the
$C^1$-regularity of the solution to the corresponding Dirichlet problem, which has been established quite recently in \cite{CFe}. 
More generally, the definition of $\DN$ via a dichotomy demands some care to adapt the different parts of the proof.

\smallskip
\item[--] Theorem \ref{t:11} gives yet another characterization of stadium-like domains, as the only convex sets  $\Omega$ where the unique solution to the homogeneous Dirichlet problem with constant source term for the normalized operator achieves its maximal regularity, namely  is in $C ^ {1, 1} (\Omega)$. 
\end{itemize}

We address as an interesting and challenging task the  problem of extending our results to non-convex domains. 

The paper is organized as follows. In Section \ref{secprel} we collect the required preliminary definitions and results. In Section \ref{secmain} we state our main results (Theorems \ref{t:s}, \ref{t:sN}, and \ref{t:11}), along with an outline of the proofs of Theorems \ref{t:s} and \ref{t:sN}, including the statement of the auxiliary results which serve as intermediate steps. In case of the operator $\Dinf$, the proofs of these intermediate steps can be found in \cite{CFd}, except for Proposition \ref{l:superweb}, which is precisely the key new ingredient allowing us to remove the diametral ball condition. In case of the operator $\DN$, the proofs of all the intermediate steps must be adapted, and therefore we have chosen to present them separately in  Section~\ref{secapp}.  Finally in Section \ref{sec11} we prove Theorem \ref{t:11}. 

\section {Preliminaries} \label{secprel}

Let us recall the basic notions and known results 
about the unique viscosity solution to the Dirichlet problems  given by the first two equations in 
\eqref{f:s} or in \eqref{f:sN}. 

For a $C^2$ function $\varphi$, we introduce the 
{\it (not normalized) infinity Laplacian} by
\[
\Dinf \varphi:= \pscal{\nabla^2\varphi\, \nabla\varphi}{\nabla\varphi}
\]
and the operators 
\[
\begin{split}
\Dp\varphi(x) & :=
\begin{cases}
%|\nabla\varphi(x)|^{-2}\, \pscal{\nabla^2\varphi(x)\, \nabla\varphi(x)}{\nabla\varphi(x)},
\DN \varphi(x),
&\text{if}\ \nabla\varphi(x)\neq 0,\\
\lmax(\nabla^2\varphi(x)),
&\text{if}\ \nabla\varphi(x)= 0,
\end{cases}
\\
\Dm\varphi(x) & :=
\begin{cases}
\DN \varphi(x),
&\text{if}\ \nabla\varphi(x)\neq 0,\\
\lmin(\nabla^2\varphi(x)),
&\text{if}\ \nabla\varphi(x)= 0 \,.
\end{cases}
\end{split}
\]
Here $\DN \varphi$ is
the {\it normalized infinity Laplacian}:
\[
\DN \varphi := \frac{1}{|\nabla \varphi |^2}\, \pscal{\nabla^2 \varphi \, \nabla \varphi}{\nabla \varphi}
\]
and, for a symmetric matrix $A\in \R ^ { n \times n }_{{\rm sym}}$, $\lmin(A)$ and $\lmax(A)$
denote respectively the minimum and the maximum eigenvalue of $A$.

Let $\Omega$ be an open bounded subset of $\R ^n$, and consider the infinity Laplace equations
\begin{equation}\label{f:infty}
-\Delta_\infty u = 1 \qquad \text{in}\ \Omega 
\end{equation}
and 
\begin{equation}
\label{f:inftyN}
-\DN u = 1\qquad \text{in}\ \Omega\,.
\end{equation}

In order to recall the notion of viscosity solutions for these equations, according to \cite{CHL},  it is convenient to fix some notation. 
If $u,v\colon\Omega\to\R$ are two functions and
$x\in\Omega$, by
\[
u \prec_x v
\]
we mean that $u(x) = v(x)$ and $u(y) \leq v(y)$ for every $y\in\Omega$.

Moreover we denote by $\osubjet u (x)$ (resp.\ $\osuperjet u (x)$) the {\it second order sub-jet} (resp.\ {\it super-jet}),
%$\osubjet u (x)$ (resp. $\osuperjet u (x)$)
of a function $u\in C(\overline{\Omega})$
at a point $x\in \overline{\Omega}$,  which is by definition the set of pairs
$(p, A) \in \R ^n \times \R ^ { n \times n }_{{\rm sym}}$ such that, as $y \to x,\ y\in \overline{\Omega}$, it holds
$$
 u (y) \geq  (\leq) \ u ( x) + \pscal{ p}{y- x} 
+ \frac{1}{2} \pscal {A (y- x)}{y- x} + o ( |y - x|^2) 
\,.
$$

A {\it viscosity solution to \eqref{f:infty}, or to  \eqref{f:inftyN}}, is a function $u\in C({\Omega})$
which is both a viscosity sub-solution and a viscosity super-solution to the same equation. 

A {\it viscosity subsolution  to  \eqref{f:infty}, or to \eqref{f:inftyN}}, is an upper semicontinuous function $u$ such that, 
for every $x \in \Omega$,  
% (meaning that $-\Delta_\infty\varphi(x)-1\leq 0$ for all $\varphi\in %C^2(\Omega)$
%such that $\varphi-u$ has a local minimum at $x$), 
$$
\forall\varphi\in C^2(\Omega)\ \text{s.t.}\ u\prec_x \varphi\,, \quad -\Delta_\infty\varphi(x)\leq 1\, , $$
or
\begin{equation}
	\label{f:subN}
	 \forall\varphi\in C^2(\Omega)\ \text{s.t.}\ u\prec_x \varphi, \ -\Dp \varphi(x) \leq 1, \text{ i.e.} \ \begin{cases}
-\Dinf \varphi(x) \leq |\nabla\varphi(x)|^2, &
\\ 
- \lmax(\nabla^2\varphi(x)) \leq 1,  \, \text{if } \nabla \varphi  (x) =0\,;&
\end{cases}
\end{equation}
equivalently, in terms of superjets, this amounts to ask respectively that
$$\forall (p, X) \in J ^ { 2 , +} _\Omega u (x) \,, \quad - \langle X p, p \rangle \leq 1\, , 
$$
or
$$
\forall(p, X) \in J ^ {2, +} _\Omega u (x)\, , \quad 
\begin{cases}
- \langle X p,  p  \rangle \leq  |p|^2, & 
\\ 
-\lambda _{\max} (X) \leq 1,  \text{ if }  p =0\,.&
\end{cases}
$$

A {\it viscosity super-solution to \eqref{f:infty}, or to \eqref{f:inftyN}},
is a lower semicontinuous function $u$ such that,  
for every $x \in \Omega$,  
$$\forall\varphi\in C^2(\Omega)\ \text{s.t.}\ \varphi\prec_x u\,, \quad -\Delta_\infty\varphi(x)\geq 1\,   $$
or
\begin{equation}
	\label{f:superN}
	\forall\varphi\in C^2(\Omega)\ \text{s.t.}\ \varphi\prec_x u, \,  -\Dm \varphi(x) \geq 1, \text{ i.e.} \
\begin{cases}
-\Dinf \varphi(x) \geq |\nabla\varphi(x)|^2, &
\\ 
- \lmin(\nabla^2\varphi(x)) \geq 1, \text { if } \nabla \varphi  (x) =0\,;&
\end{cases}
\end{equation}
equivalently, in terms of subjets, this amounts to ask that
$$
\forall (p, X) \in J ^ { 2 , -} _\Omega u (x)\, , \quad - \langle X p, p \rangle \geq 1 \, , $$
or $$
\forall(p, X) \in J ^ {2, -} _\Omega u (x)\, , \quad 
\begin{cases}
- \langle X p,  p  \rangle \geq |p|^2 &  
%\text{ if } p \neq 0
\\ 
-\lambda _{\min} (X) \geq 1\,,   \text{ if }  p =0\,.&
\end{cases}
$$

Next consider the Dirichlet boundary value problems
\begin{equation}\label{f:d}
\begin{cases} 
		-\Dinf u = 1 &\text{in}\ \Omega\,,\\
		u = 0 &\text{on}\ \partial\Omega\,, 
	\end{cases}
\end{equation}

\begin{equation}\label{f:dN}
	\begin{cases} 
		-\DN u = 1 &\text{in}\ \Omega\,,\\
		u = 0 &\text{on}\ \partial\Omega\,. 
	\end{cases}
\end{equation}

A {\it viscosity solution to} \eqref{f:d} {\it or to} \eqref{f:dN}  
is  a function $u\in C(\overline{\Omega})$ such that $u=0$ on $\partial\Omega$ and
$u$ is a viscosity solution to the pde $-\Dinf u = 1$ or $-\DN u = 1$, according to the above recalled  definitions.

The existence and uniqueness of such a viscosity solution has been proved  in \cite{LuWang, BhMo} for the Dirichlet problem \eqref{f:d}  and in \cite{PSSW, LuWang2, LuWang3, ArmSm2012} for the Dirichlet problem \eqref{f:dN}.

Concerning regularity, we proved in our previous papers \cite{CFd} and \cite{CFe} that, under the assumption that $\Omega$ is convex, the unique solution to the above Dirichlet problems is power-concave (precisely, $(\frac{3}{4})$-concave in case of problem \eqref{f:d} and $(\frac{1}{2})$-concave in case of 
problem~\eqref{f:dN}), locally semiconcave, and of class $C ^ 1 (\Omega)$. 
In case of the Dirichlet problem for the not-normalized operator, such regularity result was established 
in \cite{CFd} under the additional assumption that $\Omega$ satisfies an interior sphere condition;
we are going to remove this restriction in Lemma~\ref{l:final} below,
using the fact that an appropriate web function is a supersolution (see Proposition~\ref{l:superweb}).

Finally, we need to recall some definitions related to the distance function to the boundary of $\Omega$, which will be denoted by $d _{\partial \Omega}$. 
We let $\Sigma (\Omega)$ be the set of points in $\Omega$ where $d_{\partial \Omega}$ is not differentiable, and we call {\it cut locus} and {\it high ridge}
the sets given respectively by
\begin{eqnarray} 
\hbox{$\Cut(\Omega)$ := the closure of   $\Sigma (\Omega)$ in $\overline \Omega$
} \qquad  \qquad& \label{cut} \\ \noalign{\medskip}
\hbox{$\high (\Omega)$ := the set  where $d _{\partial \Omega}(x) = \rho _\Omega:= \max _{ \overline \Omega} d _{\partial \Omega}\,. $ } & \label{high}
\end{eqnarray}
Following \cite{gazzola, CFGa}, we say that $u:\Omega \to \R$ is a {\it web--function} if $u$
depends only on $d_{\partial \Omega}$, i.e.\
$u = g \circ d_{\partial \Omega}$ for some function $g\colon [0, \inradius] \to \R$.

Two web-functions will play a special role in the paper, in connection with 
problems~\eqref{f:s}--\eqref{f:sN}.  
We denote them by 
$ \phi ^\Omega$ and $\phi_N ^\Omega$  respectively: 
\begin{eqnarray} 
\qquad\qquad\qquad \qquad \qquad\phi ^\Omega (x) : = c_0 \left[\rho_\Omega ^{4/3} - (\rho_\Omega - \dist(x))^{4/3}\right] \,, \quad \hbox{ where } c_0 := 3^{4/3} / 4,
\,  
& \label{f:phi}
\\ \noalign{\medskip} 
\phi_N ^\Omega (x) : = \frac{1}{2} \left[\rho_\Omega ^{2} - (\rho_\Omega - \dist(x))^{2}\right]\,.\quad  \qquad \qquad \qquad\qquad \qquad
& \label{f:phiN} 
 \end{eqnarray}
%%%%%%%%%%%%%%%%%%%%%%%%%%%%%%%%%%%%%%%%%%%%%%%%%
\section{Results} \label{secmain}

{\it Throughout the paper, $\Omega$ is assumed to be an open bounded connected subset of $\R^n$. 
When the additional assumption that $\Omega$ is convex is needed, this is explicitly specified in the statements.}

\medskip

In our paper \cite{CFb}, we obtained some geometric information 
on the shape of domains $\Omega\subset \R ^n$ whose cut locus  $\Cut(\Omega)$  and high ridge $\high(\Omega)$, defined respectively in \eqref{cut} and \eqref{high}, agree. In particular we proved that, in dimension $n=2$, a domain $\Omega$ such that $\Cut(\Omega) = \high(\Omega)$ is necessarily the tubular 
neighborhood of a line segment, possibly degenerated into a point. Inspired by this characterization, we set the following 

%(Moreover, in any space dimension and without convexity assumption a domain such that $\Cut(\Omega) = \high(\Omega)$
% is necessarily a ball if $\partial \Omega$ is of class $C^2$ (see  \cite[Thm.\ 12]{CFb}). 

\begin{definition}\label{defstad} We say that an open bounded convex subset of $\R^n$ is a {\it stadium-like domain} if there holds $\Cut (\Omega) = \high (\Omega)$. 
\end{definition}

\smallskip
Our main results state that being a stadium-like domain is a necessary and sufficient condition on a convex set $\Omega$  for the existence of a solution to any of the overdetermined problems \eqref{f:s} and \eqref{f:sN}. 

\begin{theorem}\label{t:s} Assume that $\Omega$ is convex. Then  the overdetermined boundary value 
	problem~\eqref{f:s} admits a solution $u \in C ^ 1 (\overline \Omega)$ if and only if
	$\Omega$ is a stadium-like domain (and in this case  it holds $u= \phi^\Omega$, with $\rho _\Omega = c $). 
\end{theorem}

\begin{theorem}\label{t:sN} 
Assume that $\Omega$ is convex. 
Then the overdetermined boundary value problem~\eqref{f:sN} admits a solution $u \in C ^ 1 (\overline \Omega)$ if and only if 
	$\Omega$ is a stadium-like domain  (and in this case 
it holds $\uN= \phi^{\Omega}_N$, with $\rho _\Omega = c $).  
\end{theorem}

As a companion result, which will be obtained as a consequence of Theorem \ref{t:sN}, we establish that being a stadium-like domain is also a necessary and sufficient condition on a convex set $\Omega$  for the $C ^ {1, 1}$ regularity of the unique solution to the Dirichlet problem \eqref{f:dN}:

\begin{theorem}\label{t:11} 
Assume that $\Omega$ is convex. 
Then the unique solution to the Dirichlet boundary value problem~\eqref{f:dN} is of class $ C ^{ 1, 1} (\Omega)$ if and only if 
	$\Omega$ is a stadium-like domain (and in this case 
it holds $\uN= \phi^{\Omega}_N$, with $\rho _\Omega = c $).  
\end{theorem}

\begin{remark}\label{corgeo}  
By combining Theorems \ref{t:s}, \ref{t:sN} and \ref{t:11} with Theorem 6 in \cite{CFb}, we infer that,  in dimension $n=2$,  domains $\Omega$ where any of the overdetermined problems  \eqref{f:s} or \eqref{f:sN} admits a solution (or where the unique solution to problem \eqref{f:dN} is of class $C ^ {1, 1} (\Omega)$)  are geometrically characterized  as
		$$\Omega=   \{ x \in \R ^2 \ :\ {\rm dist} (x, S) < \rho_\Omega \}\,,$$
		being the set  $S:=\Cut (\Omega) = \high (\Omega)$ 
		a line segment (possibly degenerated into a point).
		If in addition $\partial \Omega$ is assumed to be of class $C ^2$, then $\Omega$ is a ball (see \cite[Theorem 12]{CFb}). 
\end{remark}

\begin{remark}\label{remreg}
The same statement as Theorem~\ref{t:11} for the not normalized operator is clearly false. In fact, notice carefully that the function $\phi ^\Omega$ is merely of class $C ^ {1, 1/3} (\Omega)$.
We recall that, in the case of infinity harmonic functions, 
the works by Savin~\cite{Sav}, Evans-Savin \cite{EvSav} and Evans-Smart \cite{EvSm}
establish they are differentiable in any space dimension and $C ^{1, \alpha}$ in dimension two.  
\end{remark}

\begin{remark}

%This geometric characterization is obtained under the assumption that the domain $\Omega$ is convex, and that the solution $u$ is $C ^ 1$ {\it up to the boundary}, meaning that 
%\begin{equation}\label{hu}
%\exists \, \delta>0 \ :\ u \text{ is of class } C ^1 \text{ on } \{ x \in \overline \Omega \ :\ {\rm dist} (x, \partial \Omega) < \delta\}\,.
%\end{equation}
We stress that asking that the solution is of class $ C ^ 1 (\overline \Omega)$  in Theorems~\ref{t:s} 
and~\ref{t:sN} amounts to require merely that the $C ^1(\Omega)$-regularity
known for the unique solution  to problems~\eqref{f:d}--\eqref{f:dN} ({\it cf.}\ Section \ref{secprel}) is preserved at $\partial \Omega$. Notice that this is somehow necessary in order to give a pointwise meaning to the Neumann boundary condition in~\eqref{f:s}--\eqref{f:sN}. 
We address as an open problem the question of establishing whether the $C^1$ regularity of the solution to problems~\eqref{f:d}-\eqref{f:dN} extends up to $\partial \Omega$  in dependence of the regularity of the boundary itself. For related boundary regularity results, see \cite{WaYu, Hong, Hong2}.
 \end{remark}

\smallskip
We now outline the  proof of Theorems \ref{t:s} and \ref{t:sN}, 
by stating the results which serve as main intermediate steps and
explaining how they allow to conclude. For convenience, the proof of such intermediate statements is postponed to Section \ref{secapp}, whereas the proof of Theorem~\ref{t:11} is given in the final Section~\ref{sec11}.

\smallskip
The main idea to prove Theorems \ref{t:s} and \ref{t:sN} is to make use of suitable $P$-functions, introduced hereafter. 

\begin{definition}\label{d:P}  For $x \in \Omega$, we set 
	\begin{equation}\label{defP}
		P (x) := \frac{|\nabla u(x)|^4}{4} + u(x) \, , \qquad \PN (x):= \frac{|\nabla \uN(x)|^2}{2} + \uN(x) \, ,	\end{equation}
where $u$ and $u _N$ denote respectively the unique solution to problems \eqref{f:d} and \eqref{f:dN}. 
\end{definition}

%\begin{definition}
%We say that $u\in C(\overline{\Omega})$ is a web--viscosity solution to~\eqref{f:dirich}
%if $u$ is a viscosity solution to~\eqref{f:dirich} that
%depends only on the distance $\dist$ from the boundary of $\Omega$, i.e.\
%$u = g \circ \dist$ for some function $g\colon [0, \inradius] \to \R$.
%\end{definition}

The choice of the above $P$-functions is due to the fact that their constancy  
on the whole $\Omega$, if satisfied, gives the crucial  information that $u$ and $u _N$ are web-functions, 
and more precisely that they agree with the functions $\phi_\Omega$ and $\phi ^N_\Omega$ introduced in \eqref{f:phi}-\eqref{f:phiN}. 
We have indeed:

\begin{proposition}
	\label{p:P1} 
	Assume that the unique viscosity solution to~problem \eqref{f:d} or~\eqref{f:dN} is of class $C ^ 1 (\Omega)$, and that
	\begin{equation}
		\label{f:pconst}
		P (x) = \lambda \ \text{a.e.\ on } \Omega \qquad \hbox{ or } \qquad \PN (x) = \lN \ \text{a.e.\ on } \Omega  \,,
	\end{equation}
where $\lambda$ and $\lN$ are positive constants satisfying $\lambda \leq c_0 \inradius^{4/3}$ and $\lN \leq  \frac{1}{2}\inradius^2 $. 
	
	Then we have respectively: $\lambda = c_0 \inradius^{4/3}$  and $u = \phi ^\Omega$, or $\lN= \frac{1}{2} \inradius^2 $ and
	$\uN=\phi _N ^\Omega$. 
\end{proposition} 

In turn, if  the unique solution to problem \eqref{f:d} or \eqref{f:dN} happens to be a web--function, we can prove that  necessarily the cut locus and high ridge of $\Omega$ agree. Actually this geometric condition turns out to be necessary and sufficient for the solution being a web--function, according to the result below: 

\begin{proposition}\label{t:websol} The unique viscosity solution to~problem \eqref{f:d} or~\eqref{f:dN}
is a web--function if and only if there holds $\Cut(\Omega) = \high(\Omega)$.
%In this case,
%the web--viscosity solution is given by
%$\phi ^ \Omega _N := g\circ\dist$, with $g$
%defined by \begin{equation}
%\label{f:g}
%g(t) := \frac{\inradius^2 - (\inradius - t)^2}{2}\,,
%\qquad t\in [0, \inradius].
%\end{equation}
\end{proposition}

In view of Propositions \ref{p:P1} and \ref{t:websol}, in order to prove Theorems \ref{t:s} and \ref{t:sN}, one is reduced to answer the following question: is it true  that, if a solution to the overdetermined 
problems~\eqref{f:s}--\eqref{f:sN} exists, the corresponding $P$-function is constant? 

In this respect, the pde interpreted pointwise at points of two-differentiability of $u$ yields an elementary but important observation. 
Let $u$ and $\uN$  be the solutions to problems~\eqref{f:d}--\eqref{f:dN}, and let
	$\gamma$ and  $\gamma _N$ be  local solutions on some interval $[0, \delta)$ to the gradient flow problems 
	\begin{equation}
	\label{f:geo}
	\begin{cases}
		\dot \gamma (t) = \nabla u (\gamma (t)) 
		%& \qquad \forall t \in [0, \delta) 
		\\
		\gamma (0) = x \in \overline \Omega\,,& 
	\end{cases}
\qquad	\begin{cases}
		\dot \gamma_N (t) = \nabla \uN (\gamma _N(t)) 
		%& \qquad \forall t \in [0, \delta) 
		\\
		\gamma _N(0) = x \in \overline \Omega\,.& 
	\end{cases}
	\end{equation}
	
We claim that, if  $u$ (resp.\ $u _N$) is twice differentiable at $\gamma (t)$ (resp.\ $\gamma _N(t)$) for $\mathcal L ^1$-a.e.\ $t \in [0, \delta)$,  then it holds
	\begin{equation}\label{f:dernulla}
		\frac{d}{dt}\big( P (\gamma (t)) \big )  = 0 \ \ \Big (\text{resp.}\ \frac{d}{dt}\big( P_N (\gamma_N (t)) \big )  = 0 \Big ) \qquad \mathcal L ^ 1\hbox{-a.e.\ in } [0, \delta) .
	\end{equation}

The proof of this claim is very simple and we limit ourselves to check it for the normalized operator, 
the other case being completely analogous. 
At every point $x$ where $u _N$ is twice differentiable, it holds
$\nabla \PN(x)= D^2\uN(x)\, \nabla \uN(x)+\nabla  \uN(x)$;
	we infer that
	\[
	\pscal{\nabla \PN(x)}{\nabla \uN (x)}
	= \Dinf \uN(x) + |\nabla \uN(x)|^2 = 0.
	\]
	Thus, since by assumption $\uN$ is twice differentiable at $\gamma_N (t)$ for $\mathcal L ^1$-a.e.\ $t \in [0, \delta)$, it holds
	\[
	\frac{d}{dt}\big( \PN (\gamma_N (t)) \big ) 
%	=\Big \langle \nabla \PN  (\gamma _N(t) ),
%	\frac{d}{dt} \gamma_N  (t) \Big\rangle 
	= \Big \langle \nabla P_N (\gamma (t)), \nabla \uN
	(\gamma _N(t) ) \Big \rangle =  0 \qquad \mathcal L ^1\text{-a.e. on } [0, \delta).
	\]
	
Unfortunately, we have not enough regularity at our disposal to infer from \eqref{f:dernulla} that the $P$-functions are constant along the gradient flows. 
In fact, since $u$ and $u _N$ need not be of class $C ^ {1, 1} (\Omega)$, the maps $t \mapsto P \circ \gamma$  and $t \mapsto P_N \circ \gamma_N$  may fail to be in $AC([0,\delta))$.  
To circumvent this lack of regularity, we argue as follows. In a first step, we proceed by finding some upper and lower bounds for the $P$-functions. They are obtained by approximating $u$ and $\uN$ by more regular functions (their supremal convolutions, see Section \ref{secapp}).

\begin{proposition}\label{t:ineqP} If $\Omega$ is  convex, 
	there holds	\begin{eqnarray}
			\min _{\partial \Omega} \frac{|\nabla u| ^ 4 }{4} \leq P (x)  \leq \max _{\overline {\Omega}} u \qquad
		\qquad \forall  x \in \overline \Omega\,, & \label{f:tesiP}
\\ \noalign{\medskip}
		\min _{\partial \Omega} \frac{|\nabla \uN| ^ 2 }{2} \leq \PN (x)  \leq \max _{\overline {\Omega}} \uN
		\qquad \forall  x \in \overline \Omega\,. & \label{f:tesiPN}
	\end{eqnarray}
\end{proposition}

The above bounds enable us to arrive at the constancy of the $P$-functions when combined with a last key ingredient, which is stated below.

\begin{proposition}\label{l:superweb} 
The
web functions $\phi ^\Omega$ and $\phi ^ \Omega _N$ are 
viscosity supersolutions respectively  to the equation $-\Dinf u = 1 $ and $- \DN u = 1$ in $\Omega$. 
\end{proposition}

\smallskip
Based on the strategy outlined above and on the preliminary results stated so far, let us give more in detail the
proof of Theorem \ref{t:sN}.  
The proof of Theorem \ref{t:s} is omitted since it is analogous, relying on the corresponding intermediate steps. The only difference is related with the removal of the interior sphere condition appearing in~\cite[Thm.\ 5]{CFd}, as mentioned in the Introduction. This is discussed in detail after the proof of Theorem \ref{t:sN}.

\smallskip
{\bf Proof of Theorem  \ref{t:sN}}.   Let $B= B _{\rho _\Omega}(x_0)$ be an inner ball of radius $\rho _\Omega$,  let $y _0$ be a fixed point in 
$ \partial B \cap \partial \Omega$ and let $\gamma$ be the line segment $[x_0, y_0]$. 
Let $\phi ^B_N$ and $\phi ^\Omega_N$ be the web--functions defined according to \eqref{f:phiN}. 
By Proposition \ref{l:superweb}, applying the comparison principle proved in ~\cite[Thm.~2.18]{ArmSm2012},
we infer that
\begin{equation}\label{f:ineq}\phi ^B_N (x) \leq \uN (x) \leq \phi ^\Omega_N (x) \qquad \forall x \in B \,.\end{equation}
We can deduce several consequences from these inequalities. 
Firstly we observe that, since both the functions $\phi ^B_N$ and $\phi ^\Omega_N$ have a relative maximum at $x_0$, by (\ref{f:ineq}) the same property holds true for $\uN$. Hence $x_0$ is a critical point of $\uN$.  In turn, we observe that 
\begin{equation}\label{f:argmax}
	x_0 \in \argmax_{\overline \Omega} (\uN )\,. 
\end{equation}
Indeed, the set of critical points of $\uN$ agrees  with the set $\argmax_{\overline \Omega} (\uN) $ where $\uN$ attains its maximum over $\overline \Omega$. This is because, by \cite[Theorem 6]{CFe}, the
	function $\uN^{1/2}$ is concave
	% (and strictly positive) 
	in $\Omega$; hence
	its gradient 
	vanishes only at maximum points of
$\uN$.

Moreover we notice that, since the distance functions $d_{\partial B}$ and $d _{\partial 
\Omega}$ agree on the line segment $\gamma$, there holds
\begin{equation}\label{f:equ} \phi ^B_N (x) = \phi ^\Omega_N (x) \qquad \forall x \in \gamma\,.\end{equation}
As a consequence of (\ref{f:ineq}) and (\ref{f:equ}), we deduce that
$\uN (x) = \phi ^\Omega_N(x) = \phi ^B_N (x)$ for all $x \in \gamma$.
Namely, there holds
\begin{equation}\label{f:udiam}
\uN (x) = 	\frac{\rho _\Omega^2 - ( \rho _\Omega - d _{\partial \Omega}(x) ) ^ 2 }{2} \qquad \forall x \in \gamma \,.
\end{equation}
It follows from (\ref{f:udiam}) that $|\nabla \uN (y_0) | =  \rho _\Omega$.  
Recalling that by assumption $\uN$ satisfies the Neumann condition $|\nabla \uN(y)| = c$ for all $y \in \partial \Omega$, 
we deduce that the value of the parameter $c$ is related to the inradius by the equality
$c =\rho _\Omega$. 
Using such equality and~\eqref{f:argmax}, we get
\[
\max_{\overline \Omega} (\uN) = \uN (x_0) = \frac{\rho _\Omega ^2}{2} = \frac{c ^ 2}{2}\,.
\]
By Proposition \ref{t:ineqP}, this implies that the $P$-function associated with $\uN$ according 
to $(\ref{defP})$ satisfies 
$$\PN (x) = \frac{c ^ 2}{2}\qquad \forall x \in\overline \Omega\,.$$ 
Since $\Omega$ is assumed to be convex, it follows 
from~\cite[Thm.\ 16]{CFe},  that  $u$ is of class $C ^ 1 (\Omega)$. Therefore, we are in a position to apply Proposition \ref{p:P1}, to obtain that $\uN = \phi^\Omega_N$ (with $\rho _\Omega = c $), and finally Proposition \ref{t:websol} to conclude that $\Cut (\Omega) = \high (\Omega)$.  \qed

\bigskip
Going through the above proof, we see that we have used all our intermediate results, stated in Propositions \ref{p:P1}, \ref{t:websol}, \ref{t:ineqP}, and \ref{l:superweb}. Since such results have been established also  in case of the not normalized operator $\Dinf$, this allows to obtain the proof of Theorem \ref{t:s}. Nevertheless, some attention must be paid, precisely when applying Proposition \ref{p:P1}, because it requires that the unique solution to problem \eqref{f:d} is of class $C ^ 1 (\Omega)$. 
Whereas in case of problem \eqref{f:dN} the $C ^ 1$-regularity of the solution  was proved in~\cite[Thm.\ 16]{CFe} for arbitrary convex domains, 
in case of problem \eqref{f:d}, it was proved in ~\cite[Cor.\ 10]{CFd} under the additional assumption that $\Omega$ satisfies an {\it interior sphere condition}.  Moreover, an inspection of the proof of~\cite[Thm.\ 16]{CFe} reveals that  it is not straightforward to adapt it  to the case of the not normalized operator.  However,  relying on the new Proposition \ref{l:superweb}, we are now  able to remove the interior sphere condition. This is done in Lemma \ref{l:final} below. It ensures that, also in case of problem \eqref{f:d}, the $C ^1$-regularity condition asked in Proposition \ref{p:P1} is fulfilled for any convex domain, thus enabling us to conclude the proof of Theorem \ref{t:s}. 

\begin{lemma}\label{l:final}
If $\Omega$ is convex, then the unique solution to problem~\eqref{f:d} is of class $C ^ 1 (\Omega)$. 
\end{lemma}

\proof 
By~\cite[Thm.\ 9]{CFd}, it is enough  to show that the unique solution to problem~\eqref{f:d} is power-concave. 
Let $u$ be such a solution. 
For $\varepsilon\in (0,1]$ let $\Omega_\varepsilon$ denote
the outer parallel body of $\Omega$ defined by
\[
\Omega_\varepsilon := \{x\in\R^n:\ \text{dist}(x, \Omega) < \varepsilon\}\,,
\]
and let $u_\varepsilon$ denote the solution to
\[
	\begin{cases} 
	-\DN u_\varepsilon = 1 &\text{in}\ \Omega_\varepsilon\,,\\
	u_\varepsilon = 0 &\text{on}\ \partial\Omega_\varepsilon\,. 
	\end{cases}
\]
Since $\Omega_\varepsilon$ satisfies an interior sphere condition
(of radius $\varepsilon$), by \cite[Cor.\ 10]{CFd}  the function
$u_\varepsilon^{1/2}$ is concave in $\Omega_\varepsilon$.
Therefore, to show that $u^{1/2}$ is concave in $\overline{\Omega}$,
it is enough to show that, as $\varepsilon \to 0$, $u_\varepsilon \to u$
uniformly in $\overline{\Omega}$.
In turn, by \cite[Thm.~5.3]{LuWang2}, this convergence
holds true provided $\left.{u_\varepsilon}\right|_{\partial\Omega}$
tends uniformly to $0$. 

To that aim we observe that, thanks to Proposition \ref{l:superweb} and the comparison principle  proved in \cite[Thm.\ 3]{LuWang}, there holds
\[
0<u _\varepsilon (x) \leq \phi^{\Omega _\varepsilon}(x) 
= c_0 \left[(\inradius + \varepsilon)^{4/3} - \inradius^{4/3}\right]
\qquad \forall x \in \partial \Omega\, , 
\]
which implies that 
$\left.{u_\varepsilon}\right|_{\partial\Omega}$
converges uniformly to $0$ on $\partial\Omega$.
\qed

% % % % % % % % % % % % % % % % % % % % % % % % % % % % % % % % % % % % % 
\section{Proofs of intermediate results}\label{secapp}

% % % % % % % % % % % % % % % % % % % % % % % % % % % % % % % % % % % % % % % % 

\subsection{Proof of Proposition \ref{p:P1}}
In case of the not normalized operator, the result has been proved in \cite[Proposition 2]{CFd}. 
Let us prove it for the normalized operator. 
%For simplicity of notation, throughout the proof we denote by $u$ in place of $\uN$ the unique solution to problem \eqref{f:dN}. 
It is clear that the constant $\lN$ is equal to $\max_{\overline \Omega} \uN$.
On the other hand, $\max_{\overline \Omega}  \uN \geq \max _{\overline \Omega}  v = \inradius^2 / 2$,
where $v$ is the radial solution of the Dirichlet problem
in a ball $B_{\rho_\Omega}\subseteq\Omega$.
Hence $\lN = \inradius^2 / 2$.

Let $H\colon\R\times\R^n\to\R$ be the Hamiltonian defined by
\[
H(u, p) := \frac{1}{2} |p|^2 + u - \lN.
\]
Then the second equality in ~\eqref{f:pconst} can be rewritten as
\begin{equation}
	\label{star}
	H(\uN(x), \nabla \uN(x)) = 0,  \qquad \mathcal L ^n \hbox{-a.e. on } \Omega\ .
\end{equation}
Since $\uN$ is of class $C^1(\Omega)$, then it follows that it is a classical
(hence also a viscosity) solution of the Dirichlet problem
\begin{equation}
	\label{f:dpH}
	\begin{cases}
		H(\uN, \nabla \uN) = 0, & \text{in}\ \Omega,\\
		\uN = 0, & \text{on}\ \partial\Omega.
	\end{cases}
\end{equation}
Since the solution to this Dirichlet problem is unique
(see {\it e.g.}\ \cite[Theorem III.1]{B}),
to prove that $\uN=\phi ^ \Omega _N$ it is enough to show that also
$\phi ^ \Omega _N$ is a viscosity solution to~\eqref{f:dpH}.
Since 
$\phi ^ \Omega _N$ is differentiable at every point
$x\in\Omega\setminus S$,
where $S := \Sigma(\Omega) \setminus \high(\Omega)$, with
\[
\nabla\phi ^ \Omega _N(x) =
\begin{cases}
(\inradius - d _{\partial \Omega}(x)) \nabla d _{\partial \Omega}(x), & \text{if}\ x\in\Omega\setminus\Sigma(\Omega),\\
0, & \text{if}\ x\in\high(\Omega),
\end{cases}
\]
we have $H(\phi ^ \Omega _N(x), \nabla\phi ^ \Omega _N(x)) = 0$
for every $x\in\Omega\setminus S$.

We remark that $\phi^\Omega_N$ is a concave function,
since it is the composition of the concave function
\begin{equation}\label{f:defg}
g(t) := \frac{1}{2} [ \inradius^2 - (\inradius - t)^2], \qquad t\in [0,\inradius],
\end{equation}
with the distance function $d _{\partial \Omega}$, which in turn is concave because
$\Omega$ is a convex set.
Since $S \subseteq \Sigma(\Omega)$ has vanishing Lebesgue measure
and $H$ is convex with respect to the gradient variable,
from Proposition~5.3.1 in \cite{CaSi} we conclude
that $\phi^\Omega_N$ is
a viscosity solution to~\eqref{f:dpH}.

\subsection{Proof of Proposition \ref{t:websol}}

For the case of problem~\eqref{f:d}, the result has been proved in 
\cite{CFc}, so we need to consider only the case of the
normalized infinity Laplacian. 
%As above, for simplicity of notation, throughout the proof we are going to denote by $u$ in place of $\uN$ the unique solution to problem \eqref{f:dN}. 

Assume that $\Omega$ is a stadium--like domain,
and let us prove that $\phi ^ \Omega _N$ is a viscosity solution to~\eqref{f:dN}.
Let $x\in\Omega$ and let us prove that	
both conditions 
\eqref{f:subN} and~\eqref{f:superN}
%in Definition~\ref{d:visc} 
are satisfied.
Let $p\in\Cut(\Omega) = \high(\Omega)$ and $q\in\partial\Omega$ be such that
$x\in [q,p]$ (we remark that we have $x=p$ if $x$ belongs itself to the
cut locus).
Let us define $\nu := (p-q) / |p-q|$.

\medskip
Let us first prove~\eqref{f:subN}.
By the comparison principle proved in \cite[Thm.~2.18]{ArmSm2012} we have 
\[
\phi ^ \Omega _N(y) \geq \frac{\inradius^2-|y-p|^2}{2} =: v(y) \,, \qquad \forall y\in \Omega,
\]
since the function $v$ at the right--hand side is the solution of the Dirichlet
problem in the ball $B_{\inradius}(p)\subseteq\Omega$.
Moreover, the functions $\phi ^ \Omega _N$ and $v$ coincide on the segment $[q,p]$
and, in particular, at the point $x$.
If $\phi ^ \Omega _N \prec_x \varphi$ we thus have
\[
\varphi(x) = \phi ^ \Omega _N(x) = v(x),
\qquad
v(y) \leq \phi ^ \Omega _N(y) \leq \varphi(y) \quad\forall y\in\Omega,
\]
so that $v \prec_x \varphi$.
Since $v$ is a solution to $-\DN v = 1$,
this implies that
$-\Dp\varphi(x)\leq 1$.
%$\lmax(\nabla^2\varphi(x)) \geq -1$.

\medskip
Let us now prove~\eqref{f:superN}.
Let $\varphi \prec_x \phi ^ \Omega _N$.

If $x\in\Cut(\Omega) = \high(\Omega)$, we must have $\nabla\varphi(x) = 0$ and
\[
\pscal{\nabla^2\varphi(x) (y-x)}{y-x} \leq - (\inradius - d _{\partial \Omega}(y))^2,
\qquad\forall y\in\overline{\Omega}.
\]
%If $y\in\partial\Omega$ is a projection of $x$, so that 
Since, in this case,
$x = q + \inradius\, \nu$
%with $\nu$ a unit vector, 
we get
\[
\lmin(\nabla^2\varphi(x))
\leq \pscal{\nabla^2\varphi(x) \nu}{\nu}
= \frac{1}{\inradius^2} \pscal{\nabla^2\varphi(x) (q-x)}{q-x} \leq -1.
\]

If $x\not\in\Cut(\Omega)$, then $\tau := d _{\partial \Omega}(x) = |x-q| < \inradius$,
and
\[
\nabla\varphi(x) = \nabla\phi ^ \Omega _N(x) = g'(\tau) \nabla\dist(x) = g'(\tau) \nu\neq 0.
\]
Moreover, we have
\[
h(t) := \varphi(q+t\nu) \leq \phi ^ \Omega _N(q + t\nu) = g(t),
\qquad\forall t\in [0, \inradius],
\]
and
\[
h(\tau) = g(\tau),\quad
h'(\tau) = g'(\tau) > 0, \quad
h''(\tau) \leq g''(\tau).
\]
In particular we get
\[
-\DN \varphi(x) = -h''(\tau) \geq -g''(\tau) = 1,
\]
so that we have proved that $\phi ^ \Omega _N$ is a super-solution to~\eqref{f:dN}.

\medskip
It remains to prove the converse implication of the proposition.
Let us assume that the unique viscosity solution to~\eqref{f:dN} is a web-function,
and let us prove that $\Omega$ is a stadium--like domain and that the solution
is given by $\phi ^ \Omega _N$.

Assume that the unique viscosity solution to~\eqref{f:dN} is of the form 
$\uN(x) := f(d_{\partial \Omega}(x))$. 

We claim that 
the map $t \mapsto f (t)$ is monotone increasing on $[0, \rho _\Omega]$,
and that the function $v(z) := f(\rho _\Omega -|z|)$ is a viscosity solution of
\begin{equation}
	\label{f:visrad}
	-\DN v = 1 \qquad \text{in}\ 
	B_{\rho_\Omega}(0)\setminus\{0\}.
\end{equation}

Namely, assume by contradiction that $t \mapsto f (t)$ is not monotone increasing  on $[0, \rho _\Omega]$: let $t_1, t_2 \in [0, \rho _\Omega]$ be such that  $t_1 < t _2$ but $f ( t _1) > f ( t_2)$. 
Then the absolute minimum of the continuous function $f$ on the interval $[t_1, \rho _\Omega]$ is attained at some point $t_0 > t _1$; in particular, there exists a point $t_0 \in (0, \rho _\Omega]$ which is of local minimum for the map $f$. 
Let us show that this  fact is not compatible with the assumption that  $\uN (x)= f ( d _{\partial \Omega} (x))$ is a web viscosity solution to $- \DN \uN = 1$ in $\Omega$.  
Since $t_0>0$, there exists a point  $x_0$ lying  in $\Omega$ such that $d _{\partial \Omega} (x_0) = t _0$. Since $t_0$ is a local minimum for the map $f$, the point $x_0$ is a local minimum for the function $\uN$. 
Then, we can construct a $C ^2$ function $\varphi$ with $\varphi \prec _{x_0} \uN$, which is locally constant in a neighbourhood of $x_0$. 
Clearly it holds $-\Dm \varphi  = - \lambda _{\min} (\nabla ^ 2 \varphi (x_0)) = 0 < 1$, against the fact that $\uN$ is a viscosity super-solution. 

To complete the proof of the claim, let us show that $v(z) := f(\rho _\Omega -|z|)$ is a viscosity solution to \eqref{f:visrad} at a fixed point $z_0\in B_{\rho_\Omega}(0)\setminus\{0\}$. If $\psi$ is
a $C^2$ function with $v \prec_{z_0} \psi$, 
we have to show that
\begin{equation}
	\label{f:sspsi}
	- \Dp\psi(z_0) 
	%=-\pscal{D^2\psi(z_0) \nabla\psi(z_0)}{\nabla\psi(z_0)} 
	\leq 1.
\end{equation}
We choose  a maximal ray $[p_0, q_0]$, with $p_0 \in M (\Omega)$ and $q_0 \in \partial \Omega$, that is, $p_0$ is the center of a ball of radius $\rho _\Omega = |p _0 - q _0|$ contained into $\Omega$. 
We pick a point $x_0  \in \Omega$ such that  
$$x_0 \in ] p _0, q _0[ \qquad \hbox{ and } \qquad d _{\partial \Omega} (x_0) = \rho _\Omega - | z_0|$$ and, for $x$ belonging to a neighborhood of $x_0$,  we set
$$
z (x) := \big [ \rho _\Omega - | x - q _0 | \big ] \zeta _0\ ,\qquad \hbox{ with } \zeta_0 := \frac{z_0} { | z_0 |}\,.
$$
In particular, notice that by construction there holds $z (x_0) = z _0$. 

We now consider the composite map
$$\varphi (x) := \psi ( z ( x))\, .$$ 
Clearly it is of class $C ^2$ in a neighborhood of $x_0$, and it is easy to check that it satisfies the condition $\uN \prec_{x_0} \varphi$.  
Indeed, by  the definitions of $\uN, v$, and $z$, and since $v \prec_{z_0} \psi$, there holds
\[
\uN(x_0) = f ( d _{\partial \Omega} (x_0) ) = f (\rho _\Omega - |z_0|) = v(z_0) = \psi(z_0) = \varphi(x_0). 
\]
Moreover there exists $r>0$ such that
\[ 
\uN(x) = f(d_{\partial \Omega}(x))  \leq f (\rho _\Omega - | z (x)| )  
= v(z(x)) \leq \psi(z(x)) = \varphi(x)
\qquad\forall x\in B_{r}(x_0). 
\]
Notice that the first inequality in the line above follows from the fact already proved that $f$ is monotone increasing, while the second one holds for $r$ sufficiently small by the assumption that $v \prec_{z_0} \psi$ and the continuity the map $z$ at $x_0$. 

Then, since $\uN \prec_{x_0} \varphi$ and by assumption $\uN$ is a viscosity solution to $- \DN \uN = 1$ in $\Omega$, 
we deduce that 
\begin{equation}
	\label{f:ssphi}
	-\Dp\varphi(x_0) 
	%= -\pscal{D^2\varphi(x_0) \nabla\varphi(x_0)}{\nabla\varphi(x_0)} 
	\leq 1.
\end{equation}
We now distinguish the two cases $\nabla \varphi (x_0) = 0$ and $\nabla \varphi (x_0)  \neq 0$. 

\medskip
\textsl {Case $\nabla \varphi (x_0)  \neq 0$.}
Setting $\delta(x) := | x - q_0|$, a direct computation yields
\[
\begin{split}
& \nabla\varphi(x) = - \pscal{\nabla\psi(z(x))}{\zeta_0}\, \nabla \delta(x),\\
& D^2\varphi(x) = \pscal{D^2\psi(z(x))\, \zeta_0}{\zeta_0}\,
\nabla \delta(x) \otimes \nabla \delta(x)
- \pscal{\nabla\psi(z(x))}{\zeta_0}\, D^2 \delta(x)\,.
\end{split}
\]
Taking into account the identities
\[
\begin{split}
& |\nabla \delta (x)| = 1, \\
& [\nabla \delta(x) \otimes \nabla \delta(x)] \nabla \delta(x) = \nabla \delta(x),\\
& D^2 \delta(x)\, \nabla \delta(x) = 0\, , 
\end{split}
\]
we obtain
\begin{equation}\label{f:1}
	\Dp \varphi(x_0) = 
	\left\langle D^2\varphi(x_0) 
	\frac{\nabla\varphi(x_0)}{|\nabla \varphi (x_0)|} ,\frac{\nabla\varphi(x_0)}{|\nabla \varphi (x_0)|} 
	\right\rangle  =
	\pscal{D^2\psi(z_0)\, \zeta_0}{\zeta_0} 
	\,.
\end{equation}
Now, from %Lemma~\ref{l:subdif} (a) 
\cite[Lemma~17(a)]{CFc}
we have
$$
\nabla\psi(z_0) = \alpha\zeta_0, \qquad \hbox{ with  } 
\alpha\in -D^+f(
\rho _\Omega -|z_0|)\,,
$$
and our current assumption $\nabla \varphi (x_0 ) \neq 0$ implies $\alpha \neq 0$. 
Therefore, \begin{equation}\label{f:2}
	\Dp \psi (z_0) = \left\langle D^2\psi(z_0) 
	\frac{\nabla\psi(z_0)}{|\nabla \psi (z_0)|} ,\frac{\nabla\psi(z_0)}{|\nabla \psi(z_0)|} \right\rangle  = 
	\pscal{D^2\psi(z_0)\, \zeta_0}{\zeta_0} \,.
\end{equation}
In view of (\ref{f:1}) and (\ref{f:2}), we conclude that, in case $\nabla \varphi (x_0) \neq 0$,  \eqref{f:sspsi} follows from \eqref{f:ssphi}.

\medskip
\textsl {Case $\nabla \varphi (x_0)  =0$.} By \eqref{f:ssphi}, we know that $- \lambda _{\max} (D ^ 2 \varphi (x_0)) \leq 1$, and we have to prove that 
$- \lambda _{\max} (D ^ 2 \psi(z_0)) \leq 1$. 
From the relation  $\nabla\varphi(x_0) = - \pscal{\nabla\psi(z_0)}{\zeta_0}\, \nabla \delta(x_0)$, we see that 
$\nabla \psi (z_0) =\alpha \zeta _0 = 0$, so that the Hessian matrices $D ^ 2 \varphi (x_0)$ and $D ^ 2 \psi (z_0)$ are related by 
\begin{equation}\label{f:relHess}
	D^2\varphi(x_0) = \pscal{D^2\psi(z_0)\, \zeta_0}{\zeta_0}\,
	\nabla \delta(x_0) \otimes \nabla \delta(x_0)\,.
\end{equation}
Since $\lambda _{\max} (D ^ 2 \varphi (x_0)) \geq -1$, there exists an eigenvector $\eta$ such that 
$\langle D ^ 2 \varphi (x_0) \eta, \eta \rangle \geq -1$.   Then \eqref{f:relHess} yields
\[
-1 \leq \langle D ^ 2 \varphi (x_0) \eta, \eta \rangle =  
\langle D ^ 2 \psi (z_0)   \zeta _0,  \zeta _0 \rangle\,
(\langle \nabla \delta (x_0), \eta \rangle )^2
\leq \langle D ^ 2 \psi (z_0)   \zeta _0,  \zeta _0 \rangle
\, , 
\]
which shows that $\lambda _{\max} (D ^ 2 \psi (z_0)) \geq -1$.

In order to prove that $v$ is a viscosity super-solution to (\ref{f:visrad}) at $z_0$, one can argue in a completely analogous way. More precisely, keeping the same definitions of $\zeta _0$, $p_0$, $q_0$, and $x_0$ as above, one has just to modify the auxiliary function $z(x)$ into
$\tilde z (x) := |x- p _0 | \zeta _0$, then replace the distance function $\delta (x)$  by $\tilde \delta (x) := |x - p _0|$, and finally
apply part (b) in place of part (a) of 
%Lemma~\ref{l:subdif}. 
\cite[Lemma~17]{CFc}.

\smallskip 
We are now ready to prove 
that $\uN$ coincides with the function $\phi^N_\Omega$ defined in~\eqref{f:phiN}.
Let $f\colon [0,\inradius]\to\R$ be a continuous function 
such that 
$\uN(x) = f(\dist(x))$. 
We have to show that  $f$ agrees with the function $g\colon [0,\inradius]\to\R$ defined by~\eqref{f:defg}.
Since $\uN$ is assumed to be a viscosity solution to the Dirichlet problem~\eqref{f:dN}, according to what proved above we know that $v(z) := f(\rho _\Omega -|z|)$ is a viscosity solution to 
\begin{equation}
	\label{f:dirisrad}
	\begin{cases}
		-\Delta^N_{\infty} v = 1 & \text{in}\ B_{\rho _\Omega}(0)\setminus\{0\},\\
		v = 0 & \text{on}\ \partial B_{\rho _\Omega}(0),\\
		v(0) = f (\rho _\Omega)\,.
	\end{cases}
\end{equation}

Let us define, for every
$r >0$, the function 
\begin{equation}\label{defgc}
	g_r(t) :=  \frac{1}{2}  \left [r ^ {2} - ( r - t ) ^ {2} \right ],
	\qquad t\in [0, r]. 
\end{equation}

We claim that there exists $r \in [\inradius, + \infty)$ such that 
\begin{equation}\label{f:rgiusto}
	g_{r} (\inradius) = f (\inradius) \,.
\end{equation}

To prove this claim, we observe that the function 
$$r \mapsto  g _ r(\inradius)= \frac{1}{2}  \left [ r ^ {2} - ( r - \rho _\Omega ) ^ {2} \right ] $$ 
maps the interval $[\inradius, + \infty)$ onto $[\frac{1}{2}  \inradius^{2}, + \infty)$. 
Thus in order to show the existence of some $ r$ such that (\ref{f:rgiusto}) holds, it is enough to prove the inequality
\begin{equation}\label{f:forigin}
	f (\inradius) \geq \frac{1}{2} \inradius^{2}\,.
\end{equation}
In turn, this inequality readily follows by the comparison principle holding for 
the Dirichlet problem~\eqref{f:dN} (see \cite[Thm.~2.18]{ArmSm2012}). 
Namely, let $x_0 \in M (\Omega)$.  
We observe that the function $w(x):= g (\rho _\Omega - |x-x_0|)$  solves $- \DN w = 1$ in $B _{\rho _\Omega}(x_0)$ and $w = 0$ on $\partial B _{\rho _\Omega}(x_0)$. 
This is readily checked since, being $w$ of class $C ^2$, it holds
$$- \DN  w (x) = 
- g'' (\inradius - |x- x_0| ) = 1 \qquad \text { for } x \neq x_0$$
and
$$\begin{cases}
- \Dp w (x_0 ) =-  \lambda _{\max} (D ^ 2 w (x_0)) = - \lambda _{\max} (-{\rm Id})  = 1\,, & \\
\noalign{\medskip}
- \Dm w (x_0 ) =-  \lambda _{\min} (D ^ 2 w (x_0)) = - \lambda _{\min} (-{\rm Id})  = 1 \,.
\end{cases}$$
On the other hand, the function $\uN$ solves $- \DN \uN = 1$ 
in $B _{\inradius}(x_0)$ and $\uN \geq 0$ on $\partial B _{\inradius}(x_0)$. 
The latter inequality can be deduced by applying
the comparison principle proved in~\cite[Thm.~2.18]{ArmSm2012}. 
Again by applying the same result, we  deduce that $\uN (x) \geq g (\rho _\Omega - |x-x_0|)$ in $B _{\rho _\Omega}(x_0)$. This implies in particular inequality \eqref{f:forigin}, as
\[
f (\rho _\Omega)  = \uN (x_0) \geq g (\rho _\Omega) = \frac{1}{2} \rho _\Omega ^ {2} \,. 
\]

Now, we have that the function
\[
g_{r}(\rho _\Omega -|z|), \qquad z\in B_{\rho _\Omega} (0),
\]
is a classical solution (and hence a viscosity solution)
to  problem \eqref{f:dirisrad}. (Notice that in particular the third equation in (\ref{f:dirisrad}) is satisfied thanks to (\ref{f:rgiusto})).

From  \cite[Thm.~2.18]{ArmSm2012}, \cite[Thm.~1.8]{LuWang2}, \cite[Cor.~1.9]{PSSW} 
%$\bem$ {\it precise statement ref} $\enm$, 
we know that
there exists a unique viscosity solution to~\eqref{f:dirisrad}.  
We conclude that, for some $r \geq \rho _\Omega$, it holds
$v(z) = g_{r} (\rho _\Omega -|z|)$, that is 
\begin{equation}\label{fgc}
	f (\rho _\Omega -|z|) = g_{r} (\rho _\Omega -|z|) \, ,
\end{equation}
or equivalently $\uN (x) = g_r ( \dist(x))$.

In order to show that $\uN = \phi ^\Omega_N$,  we are reduced to prove  that 
$r = \rho _\Omega$.   
%In turn, the latter condition holds provided $g' _r (\rho _\Omega) = 0$. 
We recall that, since $r\geq \rho_\Omega$, then
$g'_r(\rho_\Omega)\geq 0$, and that
$g'_r(\rho_\Omega) = 0$ if and only if $r = \rho _\Omega$.   
Assume by contradiction that $g' _r (\rho _\Omega)>0$.  
Let $x_0 \in \high (\Omega)$. Without loss of generality, assume that $x_0 = 0$. Thanks to the concavity of $g_r$, we have
\begin{equation}\label{f:conc1}
	\uN (x) = g _r (\dist(x)) \leq  \uN (0)+  g' _ r (\rho _\Omega)  
	(d(x) - \rho _\Omega)\,. 
\end{equation}
{}From Theorem~2 in \cite{CFd},
there exist vectors $p, \zeta \in \R ^n$, with $\pscal{\zeta}{p}\neq 0$, 
and positive constants $c, C, \delta$,
such that
\begin{equation}\label{f:estidist}
	d _{\partial \Omega}(x) \leq d _{\partial \Omega}(0) + \pscal{p}{x}
	- c \pscal{\zeta}{x}^2 +\frac{C}{2}
	|x|^2 \qquad  \forall x \in  \overline B_{\delta}(0)\,.
\end{equation}
%Let $p, \zeta \in \R ^n$ be associated with the point $x_0=0$ according to Theorem~\ref{t:estid} (applied to the %distance function $d$),
%with $\pscal{\zeta}{p}\neq 0$. 
By~\eqref{f:conc1} and~\eqref{f:estidist}, 
it holds 
\[
\uN(x) \leq \varphi(x) := \uN(0) + g_r'(\rho _\Omega)\Big [  \pscal{p}{x}
- c \pscal{\zeta}{x}^2 +\frac{1}{2 \rho _\Omega}
|x|^2
\Big ]\,,
\]
so that $\uN \prec _0 \varphi$. Since $\nabla \varphi (0) = g' _ r (\rho _\Omega) p \neq 0$, via 
some straightforward computations we obtain
\[
\Dp \varphi(0) = 
\frac{g_r'(\rho _\Omega)}{|p| ^ 2}\ \Delta_{\infty}\psi(0)
= \frac{g_r'(\rho _\Omega)}{|p| ^ 2}\ \left(
-2 c \pscal{\zeta}{p}^2 + \frac{1}{\rho _\Omega}|p|^2  \right).
\]
Since $g_r' (\rho _\Omega)>0$ and $\pscal{\zeta}{p} \neq 0$, it is enough to choose  $c>0$ large enough in order to have $\Dp\varphi(0) < -1$, contradiction.

\smallskip
Since we have just proved that $\uN = \phi ^\Omega_N$, we know that $\uN ( x) = g  ( d _{\partial \Omega}(x))$, 
with $g$ as in~\eqref{f:defg}. 
Assume by contradiction that there exists $x_0 \in \Sigma (\Omega) \setminus \high (\Omega)$. Without loss of generality, assume that $x_0 = 0$, and set $d_0 = \dist (0)$. Since we are assuming $x_0 \not \in \high (\Omega)$, it holds $d_0 < \rho _\Omega$, which implies $g' (d_0) >0$. Then, we can reach a contradiction by arguing similarly as above. Namely, thanks to the concavity of $g$, we have
\begin{equation}\label{f:conc2}
	\uN (x)  \leq \uN (0)+  g'  (d_0)  
	(\dist(x) - d_0)\,. 
\end{equation}
By~\eqref{f:conc2} and~\eqref{f:estidist}, 
we have 
\[
\uN(x) \leq \varphi(x) := \uN(0) + g'(d_0)\Big [
\pscal{p}{x}
- c \pscal{\zeta}{x}^2 +\frac{1}{2 d_0}
|x|^2
\Big ]
\, ,
\]
so that $\uN \prec _0 \varphi$. 
Moreover, since $\nabla \varphi (0) = g' (d_0) p \neq 0$, we have
\[
\Dp \varphi(0) = 
\frac{g'(d_0)}{|p| ^ 2} \ \left(
-2 c \pscal{\zeta}{p}^2 + \frac{1}{d_0}|p|^2  \right)\,.
\]
Since $g' (d_0)>0$ and $\pscal{\zeta}{p} \neq 0$, it is enough to choose  $c>0$ large enough in order to have $\Dp\varphi(0) < -1$, contradiction.  
We have thus shown that $\Sigma (\Omega ) \subseteq \high (\Omega)$. 
Since the converse inclusion holds true for all  $\Omega$, and since $M (\Omega)$ is a closed set, we conclude that the required equality $\Cut (\Omega) = \high (\Omega)$ holds.
\qed

\bigskip

\subsection{
	Proof of Proposition \ref{t:ineqP}} 
The estimates~\eqref{f:tesiP} for the not normalized infinity Laplacian
have been proved in~\cite[Thm.~4]{CFd},
so we will prove only the estimates~\eqref{f:tesiPN}
for the normalized infinity Laplacian.

To that aim we need a number of
preliminary results.
%As usual, to simplify the notation we shall denote by $u $ the unique solution
%to~\eqref{f:dN}; moreover,  
%In the sequel, 
We set for brevity
\begin{equation}\label{f:maxu}
	K := \argmax_{\overline{\Omega}} \uN\, , \qquad \mu:= \max _{\overline {\Omega}} \uN\,. 
\end{equation}

A first key step is the construction of the \textsl{gradient flow} $\X$ associated with $\uN$, and the location  
of its terminal points, according to lemma below.  
The proof is omitted since it is completely analogous to that of Lemma~3 in \cite{CFd}: we limit ourselves to mention that it is based on the local semiconcavity of $\uN$
(see \cite[Theorem 3.2 and Example 3.6]{CaYu}), which in case of the solution to problem \eqref{f:dN} has been recently  proved in \cite[Prop.~13]{CFe}.

\begin{lemma}\label{l:geo} 
	Assume that $\Omega$ is convex and that the unique solution $\uN$
	to~\eqref{f:dN} is of class $C^1(\overline{\Omega})$.
	%Under the same assumptions as in Theorem~\ref{t:sN}, 
	Then, for every point   $x\in \overline \Omega \setminus K$,	
	there exists a unique solution $\X(\cdot, x)$
	to~\eqref{f:geo}
	globally defined in $[0, +\infty)$.
	Moreover, if we set
	\begin{equation}
		\label{f:defT}
		T(x) := \sup\{t\geq 0:\ \nabla \uN (\X(t,x))\neq 0 \} \in (0, +\infty],
	\end{equation}
	then 
	\begin{equation}
		\label{f:limu}
		\lim_{t\to T(x)^-} \X(t, x)
		\in K,\qquad
		\lim_{t\to T(x)^-} \nabla \uN (\X(t,x)) = 0\,.
	\end{equation}
	Finally, there exist $x_0\in\partial\Omega$
	and $t_0\in [0, T(x_0))$ such that
	$x = \X(t_0, x_0)$.
\end{lemma}

% % % % % % % % % % % % % % % % % % % % % % % % % % % % % % % % % % % % % % % %
\medskip

As we have already mentioned in Section \ref{secmain}, the above result cannot be directly exploited to infer the constancy of $\PN$ along 
the flow $\X$, because of the possible lack of absolute continuity of $\PN$. 
In order to overcome this difficulty, we approximate $\uN$ via its supremal convolutions, defined for $\e >0$ by
\begin{equation}\label{f:ue}
	u ^ \e (x) := \sup _{y \in \R ^n} \Big \{ \tilde u (y)  - \frac { |x-y| ^ 2 }{2 \e}  \Big \} \qquad \forall x \in \R ^n\,,
\end{equation}
where $\tilde u$ is a Lipschitz extension of $\uN$ to $\R ^n$ 
with ${\rm Lip}_{\R ^n} (\tilde u) = {\rm Lip} _{\overline \Omega} (\uN)$. 

In the next lemma we state the basic properties of the functions $u ^ \e$ which we are going to use in the sequel. 
Let us recall that, 
according to \cite[Lemma 3.5.7]{CaSi}, there exists $R>0$, depending only on ${\rm Lip}_{\R ^n} (\tilde u)$, such that any point $y$ at which the supremum in (\ref{f:ue}) is attained satisfies $|y-x| < \e R$. 
Thus, setting
\begin{equation}\label{f:Ue} 
	U _\e:= \big \{ x \in \Omega \ :\ \uN (x) > \e \big \} \, , \qquad A_\e := \big \{ x \in U _\e \ :\ d_{\partial U _\e}(x) > \e R \big \}\, ,
\end{equation}
there holds
\begin{equation}\label{d:ue2}
	u ^ \e (x) = \sup _{y \in U _\e} \Big \{ \uN (y)  - \frac { |x-y| ^ 2 }{2 \e}  \Big \} \qquad \forall x \in A_\e\,. 
\end{equation}
Moreover, let us define
\begin{equation}\label{f:omegae}
	m_\e := \max_{\partial A_\e} u^\e,
	\qquad
	\Omega_\e := \{x\in A_\e : \ u^\e (x) > m_\e \}\,.
\end{equation}

\begin{lemma}\label{l:approx1} 
	Under the same assumptions of Lemma~\ref{l:geo}, 
	%on $\Omega$ as in Theorem \ref{t:sN}, 
	let $u ^\e$ and $\Omega _\e$ be defined respectively as in \eqref{f:ue} and \eqref{f:omegae}. 
	Then:
	\begin{itemize}
		\item[(i)] $u ^ \e$ is of class $C ^ {1, 1}$ on $\Omega _\e$;
		\item[(ii)] $u ^\e$ is a sub-solution to $-\DN u -1 = 0$ in $\Omega_\e$;
		\item[(iii)] as $\e \to 0 ^+$,  it holds
		\[
		\begin{array}{ll}
		&  u ^ \e \to \uN \qquad \hbox{ uniformly in } \overline \Omega, \\ \noalign{\smallskip}
		& \nabla u ^ \e \to \nabla \uN   \qquad \hbox{ uniformly in } \overline \Omega 
		\end{array}
		\]
		(so that $m_\e \to 0$ and $\Omega_\e$ converges to $\Omega$
		in Hausdorff distance). 
	\end{itemize}
\end{lemma}

\proof The proofs of  (i) and (iii) are the same as those of  the corresponding statements in \cite[Lemma~4]{CFd}.  
Let us check that also statement (ii) remains true for the normalized operator. 

Let $x \in \Omega _\e$, and let $(p, X) \in  J  ^ {2, +} _{\Omega _\e} u ^ \e (x)$. 
It follows from magical properties of supremal convolution ({\it cf.} \cite[Lemma A.5]{CHL}) that $(p, X) \in J ^ {2, +} _{\Omega _\e} \uN (y)$, where $y$ is a point at which the supremum which defines
$u ^ \e (x)$ is attained.  
Since $y \in U _\e \subset \Omega _\e$, it holds $J ^ {2, +} _{\Omega} \uN (y)= J  ^ {2, +} _{\Omega _\e} u ^ \e(x)$; therefore, we have  $(p, X) \in J ^ {2, +} _{\Omega} \uN (y)$, which implies $- \langle Xp, p \rangle\leq  1 $ in case $p \neq 0$ and $- \lambda _{\max} (X) \leq 1$ in case $p = 0$. 
\qed

\bigskip

Next we observe that, for every $\e >0$, one can consider the gradient flow $\Xe$ 
associated with $u^\e$. Namely, 
for every $x_\e\in\overline{\Omega}_\e$,
the Cauchy problem
\begin{equation}\label{f:geoe}
	\begin{cases}
		\dot \gamma_\e (t) = \nabla u ^ \e (\gamma _\e (t))\,,
		%& \qquad  \forall t \in [0, \delta_\e) 
		\\
		\gamma_\e  (0) = x _\e \in  \overline{\Omega}_{\e}\,, 
	\end{cases}
\end{equation}
admits a unique solution 
$\Xe(\cdot, x_\e)\colon [0, +\infty) \to \overline{\Omega}_{\e}$. Indeed, 
the fact that $\Xe(\cdot, x_\e)$ is
defined in $[0, +\infty)$ follows from the estimate
\[
\frac{d}{dt}u^\e (\gamma_\e(t)) 
= |\nabla u^\e(\gamma_\e(t))|^2 \geq 0,
\] 
so that $\gamma_\e (t)\in\overline{\Omega}_\e$ for
every $t\geq 0$,
while uniqueness follows from the
$C^{1,1}$ regularity of $u^\e$
stated in Lemma~\ref{l:approx1}(i).

The following lemma establishes the behavior, along the  flow $\X _\e$, 
of the approximate $P$-function defined by
\begin{equation}\label{f:Pe}
	P_\e (x) := \frac{|\nabla u^\e(x) | ^ 2}{2} + u ^\e(x) \, ,\qquad x \in \overline{\Omega}_{\e}\, ,
\end{equation}
showing that $P _\e$ increases along $\X _\e$. For the proof, we refer to \cite[Lemma 5]{CFd}.

\begin{lemma}\label{l:approx2} 
	Under the same assumptions of Lemma~\ref{l:geo},
	%as in Theorem~\ref{t:sN}, 
	let $u ^\e$, $\Omega _\e$,  and $P _\e$ be defined respectively as in \eqref{f:ue}, \eqref{f:omegae},  and \eqref{f:Pe}. 
	Then, 
	for ${\mathcal H} ^ {n-1}$-a.e.\ $x_\e \in\partial{ \Omega _\e}$,
	it holds
	\[
	P _\e(\Xe(t_1, x_\e)) \leq P _\e (\Xe (t_2, x_\e) ) \,  \qquad \forall  \, t_1, t_2  \, \hbox{ with }  0\leq t_1 \leq t_2
	\,.
	\]
	%& \Xe(\cdot, x_\e) \to \X(\cdot, x_0) 
	%\qquad \hbox{ uniformly on compact subsets of } [0, T(x_0))\,.
\end{lemma}

We are finally in a position to give the

\medskip
{\bf Proof of Proposition \ref{t:ineqP}}.  
By continuity, it is enough to show that the inequalities~\eqref{f:tesiPN} hold for all 
$x \in \Omega \setminus K$.
By Lemma~\ref{l:geo},
given $x\in\Omega\setminus K$, there exist
$x_0\in\partial\Omega $ and 
$t_0\in[0, T(x_0))$ such that $x = \X (t_0, x_0)$. 
By Lemma~\ref{l:approx2}, we may find a sequence of points
$x_\e\in \partial\Omega_\e$ converging
to $x_0$ such that, 
for every $t\geq t _0$, we have
\[
P_\e(x_\e) \leq
P_\e (\Xe(t_0, x_\e)) \leq P_\e(\Xe(t, x_\e))\,.
\] 
We now pass to the limit as $\e \to 0 ^+$ in the above inequalities: 
by using the continuous dependence for ordinary differential equations
(see e.g.\ \cite[Lemma~3.1]{Hale}), 
%it follows that
%$\Xe(\cdot, x_\e)$ converge uniformly to
%$\X(\cdot, x_0)$ on compact subsets of
%$[0, T(x_0))$. Then, 
and the uniform convergences stated in
Lemma~\ref{l:approx1}(iii), we get
\begin{equation}\label{f:pineq}
	\PN(x_0) \leq \PN (x) \leq \PN (\X(t, x_0))\,.
\end{equation}
We have
\[
\PN(x_0) = \frac{|\nabla \uN(x_0)|^2}{2}
\geq \min_{\partial\Omega} \frac{|\nabla \uN|^2}{2};
\]
on the other hand, from \eqref{f:limu}, it holds
\[
\lim_{t\to T(x_0)^-}
\PN(\X(t, x_0)) = 
\lim_{t\to T(x_0)^-}
\uN(\X(t, x_0))
\leq \mu\,.
\]
Then \eqref{f:tesiPN} follows from \eqref{f:pineq}.
\qed

\subsection{Proof of Proposition \ref{l:superweb}} The proof of this result is new for both the operators $\Dinf$ and $\DN$. Since it is analogous in the two cases, we present it just for the normalized operator. 
Let $\varphi \prec_x \phi ^ \Omega _N$.  Let $p \in M (\Omega)$ and $q \in {\partial \Omega}$ be such that $x \in [q, p]$, and let $\nu := (p-q)/|p-q|$. We distinguish three cases. 

\smallskip
{\it Case 1:} $x\in \high(\Omega)$. In this case there holds necessarily $\nabla\varphi(x) = 0$ and
\[
\pscal{\nabla^2\varphi(x) (y-x)}{y-x} \leq - (\inradius - \dist(y))^2,
\qquad\forall y\in\overline{\Omega}.
\]
%If $y\in\partial\Omega$ is a projection of $x$, so that 
Since, in this case,
$x = q + \inradius\, \nu$
%with $\nu$ a unit vector, 
we get
\[
\lmin(\nabla^2\varphi(x))
\leq \pscal{\nabla^2\varphi(x) \nu}{\nu}
= \frac{1}{\inradius^2} \pscal{\nabla^2\varphi(x) (q-x)}{q-x} \leq -1.
\]

{\it Case 2:} $x\not\in\Sigma(\Omega)$. In this case we have $\tau := d _{\partial \Omega} (x) = |x-q| < \inradius$,
and
\[
\nabla\varphi(x) = \nabla\phi ^ \Omega _N(x) = g'(\tau) \nabla\dist(x) = g'(\tau) \nu\neq 0\, , 
\]
with $g$ as in \eqref{f:defg}. 
Moreover, we have
\[
h(t) := \varphi(q+t\nu) \leq \phi ^ \Omega _N(q + t\nu) = g(t),
\qquad\forall t\in [0, \inradius],
\]
and
\[
h(\tau) = g(\tau),\quad
h'(\tau) = g'(\tau) > 0, \quad
h''(\tau) \leq g''(\tau).
\]
In particular we get
\[
-\DN \varphi(x) = -h''(\tau) \geq -g''(\tau) = 1. 
\]
{\it Case 3}: $x \in \Sigma(\Omega) \setminus  M (\Omega)$. In this case, since the sub-differential of $d _{\partial \Omega}$ at $x$ is empty, the same holds true for the sub-differential of $\phi ^\Omega_N $ at $x$.  In particular, the second order sub-jet $J ^ {2, -} _\Omega \phi ^\Omega _N(x)$ is empty, so that $\phi ^\Omega_N$ trivially satisfies the definition of viscosity super-solution at $x$. 
\qed
%%%%%%%%%%%%%%%%%%%%%%%%%%%%%%%%%%
\section{Proof of Theorem \ref{t:11}}\label{sec11}

% % % % % % % % % % % % % % % % % % % % % % % % % %

The {\it sufficiency} part in the statement of Theorem \ref{t:11} readily follows from Theorem \ref{t:sN} and formula \eqref{f:phiN}. 
The {\it necessary} part is proved in Proposition \ref{p:notreg1} below, after the following preliminary lemma.

% % % % % % % % % % % % % % % % % % % % %
\begin{lemma}\label{l:geo2}

Assume that $\Omega$ is convex and that the unique solution $u$ to problem \eqref{f:dN} belongs to $C^{1,1}({\Omega}\setminus K)$, with $K $ as in \eqref{f:maxu}.  
Then, for a.e.\ $x \in \Omega \setminus K$, 
there
exists a unique solution $\X(\cdot, x)$ to \eqref{f:geo}, globally defined in $[0, + \infty)$, which satisfies 
\begin{equation}\label{f:limu2}
		 \X(t,x) \not \in K  \quad \forall t \in [0, + \infty) 
	\end{equation}
	and
\begin{equation}\label{f:limu3}	
\lim _ { t \to + \infty} {\rm dist} (\X (t, x) , K ) = 0\,.
\end{equation}
\end{lemma}

\proof
For every $x\in\Omega\setminus K$,
any local solution $\gamma$ to the second Cauchy problem in \eqref{f:geo} cannot
exit from $\{u \geq u(x)\}$ because we have
\begin{equation}
\label{f:monot}
\frac{d}{dt} u(\gamma(t))
= \nabla u(\gamma(t))\cdot\dot{\gamma}(t)
= |\nabla u(\gamma(t))|^2  \, ,  
\end{equation}
so that $u$ increases along the flow. 
Hence local solutions are actually global solutions, i.e.\ they are defined on $[0, + \infty)$. 
% NOTA: non � detto che $u\in C^1$, ma l'esistenza globale
% e' garantita dal fatto che $\nabla u = 0$ su $K$ 
The uniqueness of the gradient flow associated with $u$
in ${\Omega}\setminus K$ follows from the local
Lipschitz regularity of $\nabla u$ assumed therein. 

Let us now prove that, for a.e.\ $x \in \Omega\setminus K$, 
condition \eqref{f:limu2} is fulfilled. 
To that aim we are going to exploit the following claim, where the constant $\mu$ is defined according to \eqref{f:maxu}:

\smallskip \textsl {Claim: 
There exists a  set $L \subseteq (0, \mu)$ with $|L| = \mu$ such that, for all $m \in L$, condition \eqref{f:limu2} is satisfied for 
$\mathcal H ^ {n-1} \hbox{-a.e.}\ x \in \{ u = m \}$.} 

\smallskip
Let us first show how the lemma follows from the claim. We point out that  the set $F$ given by points $x \in \Omega \setminus K$ such that \eqref{f:limu2} is {\it false} is $\mathcal L^n$-measurable. Indeed, $F$ is open because  its complement is given by  $\bigcap _n G _n$, with 
$$G_n  := \Big \{ x \in \Omega \setminus K \ : \ \X ( t , x ) \not\in K \ \ \forall t \in [ 0 , n ] \Big \}\,,$$ 
and every $G_n$ is closed by continuous dependence on initial data.
Then, we can integrate $|\nabla u|$ over $F$ and  we obtain
\begin{equation}\label{f:coarea2}
 \int _{F}  |\nabla u| \, dx = \int _0 ^ {\mu} \, dm \int _{\{ u  = m \} \cap F } \, d \mathcal H ^ {n-1} (y) =  0\,,
\end{equation}
where the first equality holds by the coarea formula, and the second one is consequence of our claim. 

We now observe that  $|\nabla u|>  0$ 
on $\Omega \setminus K$: this is due to the fact that $u \in C ^ 1 (\Omega)$ with $u ^ {1/2}$ concave, so that $\nabla u$ vanishes only at maximum points of $u$. 

In view of this observation, \eqref{f:coarea2} implies  that $F$ is $\mathcal L^n$-negligible, and the lemma is proved. 

\smallskip 

Finally, let us give the

\smallskip
\textsl {Proof of the Claim:} 
Let us define $L$ as the set of values $m \in (0, \mu)$ such that  $u$ is twice differentiable $\mathcal H ^ {n-1}\hbox{a.e.}$ on $\{ u = m\}$.

Firstly let us check that  $|L| = \mu$. Namely, by the coarea formula, if $Z$ is the set of points in $\Omega \setminus K$ where $u$ is not twice differentiable,   we have 
$$
0 = \int _{Z}  |\nabla u| \, dx = \int _0 ^ {\mu} \, dm \int _{\{ u  = m \} \cap Z } \, d \mathcal H ^ {n-1} (y) \,.
$$
We infer that, for ${\mathcal L}^1$-a.e. $m \in (0, \mu)$, the set $\{ u  = m \}\cap Z$ is $\mathcal H ^ {n-1}$-negligible, so that $L$ is of full measure in $(0, \mu)$.

From now on, let $m$ denote a fixed value in $L$. 
For $x \in \{ u = m \}$, let us define
\begin{equation}\label{f:N}
N (x) :=  \Big \{ t \in [0, T(x)] \ :\ u   \text{ is not twice differentiable at }  \X (t, x) \Big \}
\end{equation}
and let us show that
\begin{equation}\label{f:NN}
\mathcal{L}^1(N(x)) = 0 \qquad \hbox{ for }
\mathcal{H}^{n-1}\hbox{-a.e.}\ x \in \{ u = m \}\,.
\end{equation}

By construction the set
\[
E:= \Big \{ \X (t, x) \ :\ x \in  \{ u = m \} \, ,\ t \in N (x) \Big \}
\]
is contained into the set of points where $u$ is not twice differentiable.
Then, since by assumption $u \in C ^ {1, 1} (\Omega \setminus K)$, the set $E$ is Lebesgue negligible. By the area formula, we have
\[
0 = {\mathcal L } ^n (E) = \int _{\{ u = m \}} \, d {\mathcal H} ^ {n-1} (x) 
\int _{N (x)} J \X(t, x) \, dt  \,,
\]
where $J \X$ is the Jacobian of the function $\X$
with respect to the second variable.
Since this Jacobian is strictly positive ({\it cf.}\ \cite[eq.\ (5)]{AmCr}), we infer that (\ref{f:NN}) holds true.

Let us prove that \eqref{f:limu2} holds
for every $x_0\in \{ u = m \}$ such that both the conditions
$\mathcal{L}^1(N(x_0)) = 0$ and $u$ twice differentiable at $x_0$ hold. 

Let $x_0$ be such a point, and set
%\begin{equation}
%		\label{f:defT2}
%		T(x_0) := \sup\{t\geq 0:\ \X(t,x_0) \not \in K  \} \, . 
%	
%Since $u(x_0 ) = m < \mu$, it cannot be $\nabla u (x_0) = 0 $, so that 
%$T (x_0 )$ is strictly positive. 
%
%
%Let us show that $T (x_0 ) = + \infty$. 
\[
\varphi(t) := u(\X(t, x_0)) \,,\qquad t \in [0 + \infty)\,.
\]
Since $\mathcal{L}^1(N(x_0)) = 0$, and since $u$ is assumed to be in $C ^ {1, 1} (\Omega \setminus K)$, the $P$-function is constant along $\gamma$. Therefore, 
the function
$\varphi (t)$ 
(which is in $AC ([0, + \infty))$, because
$u\in C ^ 1 (\Omega)$ 
 and $\gamma\in AC ([0, + \infty))$),
solves the Cauchy problem
\[
\begin{cases}
\frac {d \varphi}{dt} (t) = 2\lambda - 2\varphi(t) \qquad \mathcal L ^1\text{-a.e. on } [0, + \infty)\\
\varphi(0) = m\, , 
\end{cases}
\]
where $m := u (x_0)$. Since this Cauchy problem admits
a unique global solution, given  by 
$$\overline\varphi(t) := m e^{-2t} + \lambda (1 - e^{-2t})\,,$$ 
we conclude that $u(\gamma(\cdot))$ agrees with $\overline\varphi (\cdot)$ on $[0, + \infty)$. 

We now observe that 
$$\frac {d \overline \varphi (t)}{dt}   \neq 0  
\qquad \forall t \in [0, + \infty)\,.$$ 
Since $\nabla u = 0$ on $K$, we infer that
 $\X(t, x_0)\not\in K$ for $t\in [0, + \infty)$.

Eventually, we observe that \eqref{f:limu2} implies \eqref{f:limu3}. Namely,  
assume that \eqref{f:limu3} is false.  Since $u$ is increasing along the flow, there exists some level set $\{ u \leq m \}$, with $m < \mu$,  which contains the whole trajectory $\X (t, x)$ for $t \in [0, + \infty)$. 
On the compact set $\{u \leq m \}$, the continuous function $|\nabla u|$ is bounded below by some strictly positive constant. Then, in view of \eqref{f:monot}, we deduce that \eqref{f:limu2} cannot hold. 
 \qed

\bigskip

\begin{proposition}\label{p:notreg1}
Assume that $\Omega$ is convex.  If
the unique solution to problem \eqref{f:dN} is in $C^ {1,1} (\Omega \setminus K)$, then $\Omega$ is a stadium-like domain. 
\end{proposition}

\proof Let $u$ denote the unique solution to problem \eqref{f:dN}. 
As a first step we observe that, since by assumption $u \in C^ {1,1} (\Omega \setminus K)$, there holds 
\begin{equation}\label{f:Pconst}
\PN (x) = \mu \qquad \forall x \in \Omega\,.
\end{equation}

This can be obtained as a consequence of Lemma \ref{l:geo2}, by arguing as follows.  
Since $\PN$ is continuous in $\Omega$, it is enough to prove that the equality $\PN (x) = \mu$ holds a.e.\ on $\Omega \setminus K$. Namely, let us show that it holds for every $x \in \Omega \setminus K$ such that \eqref{f:limu2}-\eqref{f:limu3} hold and $\mathcal L ^ 1 (N(x)) = 0$, with $N(x)$ as in \eqref{f:N}. (Actually, both these conditions are satisfied up to a $\mathcal L^n$-negligible set, by the same arguments used in the proof of Lemma \ref{l:geo2}). Let $x \in \Omega \setminus K$ be such that \eqref{f:limu2}-\eqref{f:limu3} hold and $\mathcal L ^ 1 (N(x)) = 0$. 
Since $\mathcal L ^ 1 (N(x)) = 0$, $P$ is contant along $\gamma$ and, since \eqref{f:limu2}-\eqref{f:limu3} hold, the constant is precisely equal to $\mu$, yielding \eqref{f:Pconst}. 

Now, for $m>0$, consider the (convex) level sets
$$\Omega _m := \Big \{ x \in \Omega \ :\ u (x ) >m \Big \}\,.$$
As a consequence of \eqref{f:Pconst}, and since $u \in C ^ 1 (\overline {\Omega _m})$, $u$ satisfies on $\Omega_m$ the overdetermined boundary value problem
\begin{equation}\label{f:Sm}
\begin{cases} 
		-\DN u = 1 &\text{in}\ \Omega_m,\\
		u = m &\text{on}\ \partial\Omega_m,\\ 
		|\nabla u| = \sqrt {2(\mu -m)} &\text{on}\ \partial\Omega_m\,. 
	\end{cases}
\end{equation}
By applying Theorem \ref{t:sN} (to the function $u-m$), we infer that $\Omega _m$ is a stadium-like domain  for every $m>0$. 

To conclude, we notice that  $\{\Omega _m\}$ is an increasing sequence of open sets contained into a fixed ball; therefore,  as $m \to 0^+$, it converge in Hausdorff distance to their union (see for instance 
\cite[Section 2.2.3]{HP}). Taking into account that $\Omega = \{ u >0 \} = \bigcup_m \Omega _m$, we infer that $d_H (\Omega _m, \Omega) \to 0$, so that also the limit set $\Omega$ is a stadium-like domain. \qed

\bigskip
\textsc{Acknowledgments.}
The authors have been supported by the Gruppo Nazionale per l'Analisi Matematica, 
la Probabilit\`a e le loro Applicazioni (GNAMPA) of the Istituto Nazionale di Alta Matematica (INdAM).

%% % % % % % % % % % % % % % % % % % % % % % % % % % % % % % % % % % % % % % %
%
%\bibliographystyle{mybst}
%\bibliography{Ricerca,Graziano}

\begin{thebibliography}{10}
	
	\bibitem{AmCr}
	{L.} Ambrosio and {G.} Crippa, \emph{Existence, uniqueness, stability and
		differentiability properties of the flow associated to weakly differentiable
		vector fields}, Transport equations and multi-{D} hyperbolic conservation
	laws, Lect. Notes Unione Mat. Ital., vol.~5, Springer, Berlin, 2008,
	pp.~3--57.
	
	\bibitem{ArmSm2012}
	{S.N.} Armstrong and {C.K.} Smart, \emph{A finite difference approach to the
		infinity {L}aplace equation and tug-of-war games}, Trans.\ Amer.\ Math.\ Soc.
	\textbf{364} (2012), no.~2, 595--636. \MR{2846345}
	
	\bibitem{B}
	G.~Barles, \emph{Uniqueness and regularity results for first-order
		{H}amilton-{J}acobi equations}, Indiana Univ. Math. J. \textbf{39} (1990),
	no.~2, 443--466.
	
	\bibitem{BhMo}
	{T.} Bhattacharya and {A.} Mohammed, \emph{Inhomogeneous {D}irichlet problems
		involving the infinity-{L}aplacian}, Adv. Differential Equations \textbf{17}
	(2012), no.~3-4, 225--266.
	
	\bibitem{BH}
	F.~Brock and A.~Henrot, \emph{A symmetry result for an overdetermined elliptic
		problem using continuous rearrangement and domain derivative}, Rend. Circ.
	Mat. Palermo (2) \textbf{51} (2002), no.~3, 375--390.
	
	\bibitem{BrPr}
	{F.} Brock and {J.} Prajapat, \emph{Some new symmetry results for elliptic
		problems on the sphere and in Euclidean space}, Rend.\ Circ.\ Mat.\ Palermo
	\textbf{49} (2000), 445--462.
	
	\bibitem{butkaw}
	{G.} Buttazzo and {B.} Kawohl, \emph{Overdetermined boundary value problems for
		the {$\infty$}-{L}aplacian}, Int. Math. Res. Not. IMRN (2011), 237--247.
	
	\bibitem{CaSi}
	{P.} Cannarsa and {C.} Sinestrari, Semiconcave functions, {H}amilton-{J}acobi
	equations and optimal control, Progress in Nonlinear Differential Equations
	and their Applications, vol.~58, Birkh\"auser, Boston, 2004.
	
	\bibitem{CaYu}
	{P.} Cannarsa and {Y.} Yu, \emph{Singular dynamics for semiconcave functions},
	J. Eur. Math. Soc. (JEMS) \textbf{11} (2009), no.~5, 999--1024.
	
	\bibitem{CHL}
	{M.G.} Crandall, {H.} Ishii, and {P.L.} Lions, \emph{User's guide to viscosity
		solutions of second order partial differential equations}, Bull. Amer. Math.
	Soc. (N.S.) \textbf{27} (1992), 1--67.
	
	\bibitem{CFe}
	{G.} Crasta and {I.} Fragal\`a, \emph{A $C ^1$ regularity result for the
		inhomogeneous normalized infinity Laplacian}, to appear in Proc.\ Amer.\
	Math.\ Soc., 2015.
	
	\bibitem{CFb}
	{G.} Crasta and {I.} Fragal\`a, \emph{On the characterization of some classes
		of proximally smooth sets}, to appear in ESAIM Control Optim.\ Calc.\ Var.,
	2015.
	
	\bibitem{CFd}
	{G.} Crasta and {I.} Fragal\`a, \emph{On the Dirichlet and Serrin problems for
		the inhomogeneous infinity Laplacian in convex domains: Regularity and
		geometric results}, to appear in Arch.\ Rational Mech.\ Anal., 2015.
	
	\bibitem{CFc}
	{G.} Crasta and {I.} Fragal\`a, \emph{A symmetry problem for the infinity
		Laplacian}, to appear in Int. Mat. Res. Not. IMRN, 2015.
	
	\bibitem{CFGa}
	{G.} Crasta, {I.} Fragal\`a, and {F.} Gazzola, \emph{A sharp upper bound for
		the torsional rigidity of rods by means of web functions}, Arch.\ Rational
	Mech.\ Anal. \textbf{164} (2002), 189--211.
	
	\bibitem{EvSav}
	{L.C.} Evans and {O.} Savin, \emph{{$C^{1,\alpha}$} regularity for infinity
		harmonic functions in two dimensions}, Calc. Var. Partial Differential
	Equations \textbf{32} (2008), 325--347.
	
	\bibitem{EvSm}
	{L.C.} Evans and {C.K.} Smart, \emph{Everywhere differentiability of infinity
		harmonic functions}, Calc. Var. Partial Differential Equations \textbf{42}
	(2011), 289--299.
	
	\bibitem{f}
	{I.} Fragal{\`a}, \emph{Symmetry results for overdetermined problems on convex
		domains via {B}runn-{M}inkowski inequalities}, J. Math. Pures Appl. (9)
	\textbf{97} (2012), no.~1, 55--65.
	
	\bibitem{fg}
	{I.} Fragal{\`a} and {F.} Gazzola, \emph{Partially overdetermined elliptic
		boundary value problems}, J. Differential Equations \textbf{245} (2008),
	1299--1322.
	
	\bibitem{fgk}
	{I.} Fragal{\`a}, {F.} Gazzola, and {B.} Kawohl, \emph{Overdetermined problems
		with possibly degenerate ellipticity, a geometric approach}, Math. Z.
	\textbf{254} (2006), 117--132.
	
	\bibitem{GL}
	{N.} Garofalo and {J.L.} Lewis, \emph{A symmetry result related to some
		overdetermined boundary value problems}, Amer. J. Math. \textbf{111} (1989),
	9--33.
	
	\bibitem{gazzola}
	{F.} Gazzola, \emph{Existence of minima for nonconvex functionals in spaces of
		functions depending on the distance from the boundary}, Arch. Ration. Mech.
	Anal. \textbf{150} (1999), 57--75.
	
	\bibitem{Hale}
	{J.K.} Hale, Ordinary differential equations, second ed., Robert E. Krieger
	Publishing Co. Inc., Huntington, N.Y., 1980.
	
	\bibitem{HP}
	{A.} Henrot and {M.} Pierre, Variation et Optimisation de Formes. Une Analyse
	G\'eom\'etrique, Springer, Berlin, 2005.
	
	\bibitem{Hong}
	{G.} Hong, \emph{Boundary differentiability for inhomogeneous infinity
		{L}aplace equations}, Electron. J. Differential Equations (2014), No. 72, 6.
	
	\bibitem{Hong2}
	{G.} Hong, \emph{Counterexample to {$C^1$} boundary regularity of infinity
		harmonic functions}, Nonlinear Anal. \textbf{104} (2014), 120--123.
	
	\bibitem{K2}
	B.~Kawohl, \emph{Overdetermined problems and the {$p$}-{L}aplacian}, Acta Math.
	Univ. Comenian. (N.S.) \textbf{76} (2007), 77--83.
	
	\bibitem{KoSe}
	{R. V.} Kohn and {S.} Serfaty, \emph{A deterministic-control-based approach to
		motion by curvature}, Comm. Pure Appl. Math. \textbf{59} (2006), no.~3,
	344--407.
	
	\bibitem{LuWang}
	{G.} Lu and {P.} Wang, \emph{Inhomogeneous infinity {L}aplace equation}, Adv.
	Math. \textbf{217} (2008), 1838--1868.
	
	\bibitem{LuWang2}
	{G.} Lu and {P.} Wang, \emph{A {PDE} perspective of the normalized infinity
		{L}aplacian}, Comm. Partial Differential Equations \textbf{33} (2008),
	no.~10-12, 1788--1817. \MR{2475319 (2009m:35150)}
	
	\bibitem{LuWang3}
	{G.} Lu and {P.} Wang, \emph{Infinity {L}aplace equation with non-trivial
		right-hand side}, Electron. J. Differential Equations (2010), No. 77, 12.
	
	\bibitem{PSSW}
	{Y.} Peres, {O.} Schramm, {S.} Sheffield, and {D. B.} Wilson, \emph{Tug-of-war
		and the infinity {L}aplacian}, J. Amer. Math. Soc. \textbf{22} (2009), no.~1,
	167--210.
	
	\bibitem{Sav}
	O.~Savin, \emph{{$C^1$} regularity for infinity harmonic functions in two
		dimensions}, Arch. Ration. Mech. Anal. \textbf{176} (2005), no.~3, 351--361.
	\MR{2185662 (2006i:35108)}
	
	\bibitem{Se}
	{J.} Serrin, \emph{A symmetry problem in potential theory}, Arch.\ Rational
	Mech.\ Anal. \textbf{43} (1971), 304--318.
	
	\bibitem{Vog}
	{A.L.} Vogel, \emph{Symmetry and regularity for general regions having a
		solution to certain overdetermined boundary value problems}, Atti Sem. Mat.
	Fis. Univ. Modena \textbf{40} (1992), 443--484.
	
	\bibitem{WaYu}
	{C.} Wang and {Y.} Yu, \emph{{$C^1$}-boundary regularity of planar infinity
		harmonic functions}, Math. Res. Lett. \textbf{19} (2012), no.~4, 823--835.
	
\end{thebibliography}
%
%
%
%\end{document}
%

\def\cprime{$'$}
\providecommand{\bysame}{\leavevmode\hbox to3em{\hrulefill}\thinspace}
\providecommand{\MR}{\relax\ifhmode\unskip\space\fi MR }
\providecommand{\MRhref}[2]{%
	\href{http://www.ams.org/mathscinet-getitem?mr=#1}{#2}
}
\providecommand{\href}[2]{#2}

\end{document}